\documentclass[12pt]{article}
\usepackage{amsmath,amsfonts,amssymb,amsthm}

\def\noi{\noindent}

\def\pf{\noi{\bf Proof.\ \,}}
\def\eop{{$\square$}}

\def\labtt#1{\label {#1}}
\def\labttr#1{\label {#1}\rm }

\def\reftt#1{\ref{#1}}

\def\a{\alpha}

\def\l{\lambda}

\def\vep{\varepsilon}

\def\CC{{\mathbb C}}
\def\FF{{\mathbb F}}
\def\QQ{{\mathbb Q}}

\def\ZZ{{\mathbb Z}}

\def\la{\langle}
\def\ra{\rangle}
\def\<{\langle}
\def\>{\rangle}

\def\bs{\it}            
\def\Aut{{\bs Aut}}

\def\l{{\lambda}}

\def\half{{1 \over 2}}
\def\fourth{{1 \over 4}}

\def\dual#1{#1^*}        

\def\kron#1#2{\delta_{#1#2}}  

\def\leh{L_{E_8}}
\def\rtleh{\sqrt 2 L_{E_8}}

\def\lmt{L^-(t)}

\def\lpt{L^+(t)}

\def\lvept{L^\vep (t)}

\def\weh{W_{E_8}}

\def\ratholoex#1{2^{1+2#1}_+\Omega^+(2#1,2)}

\def\dg#1{{\cal D}({#1})}  

\def\bw#1{BW_{{2^{#1}}}}

\def\brw#1{BRW^+(2^{#1})}

\def\gd#1{G_{2^{#1}}} 
\def\rd#1{R_{2^{#1}}}

\begin{document}

\newtheorem{thm}{Theorem}[section]
\newtheorem{prop}[thm]{Proposition}
\newtheorem{lem}[thm]{Lemma}
\newtheorem{rem}[thm]{Remark}
\newtheorem{coro}[thm]{Corollary}
\newtheorem{conj}[thm]{Conjecture}
\newtheorem{de}[thm]{Definition}
\newtheorem{hyp}[thm]{Hypothesis}

\newtheorem{nota}[thm]{Notation}
\newtheorem{ex}[thm]{Example}
\newtheorem{proc}[thm]{Procedure}

\centerline{\Large \bf   Involutions on the the Barnes-Wall lattices  }

\centerline{\Large \bf   and their fixed point sublattices, I.   }

\vskip 1cm 
  
\centerline{version 21 July, 2005}
\begin{center}

\vspace{10mm}
Robert L.~Griess Jr.
\\[0pt]
Department of Mathematics\\[0pt] University of Michigan\\[0pt]
Ann Arbor, MI 48109 USA \\[0pt]

\end{center}

\vskip 0.5cm 
\begin{abstract}  We study the sublattices of the rank $2^d$  Barnes-Wall lattices $\bw d$ which occur as fixed points of involutions.  They have ranks  $2^{d-1}$ (for dirty involutions) or 
$2^{d-1}\pm 2^{k-1}$ (for clean involutions), where 
$k$, the defect,  is an integer at most $\frac d 2$.  
We discuss the involutions on $\bw d$ and determine the isometry groups of the fixed point sublattices for all involutions of defect 1.  
Transitivity results for the Bolt-Room-Wall group on isometry types of  sublattices  extend those in \cite{bwy}.  
Along the way, we classify the orbits of $AGL(d,2)$ on the Reed-Muller codes $RM(2,d)$ and describe  {\it cubi sequences} for short codewords, which give them as  Boolean sums of codimension 2  affine subspaces.   
\end{abstract}

\tableofcontents


\section{Introduction}

We continue to study the Barnes-Wall lattices $\bw d$ and their isometry groups, which are the Bolt-Room-Wall groups $\brw d \cong \ratholoex d$ for $d\ge 2, d \ne 3$ and $\weh$ for $d=3$.  In particular, we classify involutions in $\brw d$ and determine properties of their fixed point sublattices, including automorphism groups.  
For background, we analyze words of the Reed-Muller code $RM(d,2)$ in some detail and in particular determine the orbits of $AGL(d,2)$.  

We shall be using the Barnes-Wall-Ypsilanti uniqueness theory as developed in \cite {bwy}.   We recommend this article for background and terminology.  
{\it Notational warning:}  $O(L)$ means orthogonal group on a quadratic space $L$ but $O(G)$ means $O_{2'}(G)$ for a finite group $G$.

The main results of this article are described below.  See \reftt{rm}, \reftt{midset} 

\begin{thm}\labtt{maintheorem1}  The orbits for the action of $AGL(d,2)$ on the Reed-Muller code $RM(2,d)$ are as follows (for each category, there is one orbit for each allowed value of $k$): 

Short sets of defect $k=0, \dots , \lfloor \frac {d} 2 \rfloor$, which are of the form $S_1+\dots +S_k$,  where the $S_i$ are affine codimension 2 spaces which are linearly coindependent with respect to an origin in their common intersection; such a set has cardinality (or Hamming weight) $2^{d-1}-2^{d-k-1}$.  

Long sets, which are complements of short sets.

Midsets, of cardinality $2^{d-1}$, which are either affine hyperplanes (defect $0$) or nonaffine midsets of the form $S+H$, where $H$ is an affine hyperplane and $S$ is a short set of weight $2^{d-1}-2^{d-k-1}$, for a unique $k \in \{ 1, \dots , \lfloor \frac {d-1} 2 \rfloor \}$.    (Note: $k\ne \frac d 2$ here.)  
\end{thm}

Some background in the structure of BRW groups is required to state our main results.  We refer the reader to the Appendix for a summary and notations.
For definitions of clean and dirty, see \reftt{cleandirty} and for defect, see \reftt{defect}.

\begin{thm}\labtt{maintheorem2}  (i) When $d$ is odd, 
the conjugacy classes for involutions in the BRW group $\brw d$ are 
represented by the transformations:  

(Split Case) $\vep_X$, where $X$ is a codeword as listed in \ref{maintheorem1}, one for each value of the defect, $k\le \frac {d-1}2$.

(Nonsplit Case) $\eta_{d,2k,\vep}$, for $k=1, \dots ,  \frac {d-1}  2 $, $\vep=\pm$.  

(ii) 
When $d$ is even, 
the conjugacy classes for involutions in the BRW group $\brw d$ are 
represented by the transformations:  

(Split Case) $\vep_X$, where $X$ ranges over the codewords listed in \ref{maintheorem1}, but one for each value of the defect, $k$, 
together with the single clean involution $\vep_Y^\tau$, where $Y$ is a short codeword with defect $k=\frac d 2$ and $\tau$ is an outer automorphism of $\brw d$.  

(Nonsplit Case) $\eta_{d,2k,\vep}$, for $k=1, \dots ,  \frac d 2 $, where $\vep=\pm$ except for $k=\frac d 2$ when $\vep=+$ only.   
\end{thm}

The next result extends transitivity results in \cite{bwy} to a wider class of sublattices.

\begin{proc}\labttr{conjfpsublattices}
{\bf (Conjugacy for involution fixed point sublattices and recognition criteria for such.)}   
Two RSSD sublattices $M_1, M_2$  of $\bw d$ are in the same orbit of 
$\gd d$ if and only if their associated involutions  are conjugate.   We may use \reftt{maintheorem2} as a guide to orbits of $\brw d$ on RSSD sublatttices.  In particular, whether two given RSSD sublattices are in the same orbit of $\brw d$ may be decided {\it within the lattice}
by surveying a family of RSSD sublattices of $\bw d$.  It is unnecessary to examine the explicit representation of the group $\brw d$.  
See \ref{distinguishdirty}.  
\end{proc}

\begin{de}\labttr{inh} 
In general, if $X$ is a subobject of $Y$, the {\it inherited group} means the image in $Sym(X)$ of $Stab_{Aut(Y)}(X)$.  
\end{de}  

In the next result, this applies to the containment $\lvept \le L:=\bw d$.  

\begin{thm}\labtt{maintheorem4}
Consider a clean involution $t$  of defect 1 on $L:=\brw d$.

When the trace of $t$ is positive, the rank of $\lpt$ is $2^{d-2}3$.  The automorphism group is inherited when $d\ge 2, d\ne 3$ and for $d=3$ it is 
$W_{B_6}$.  

When the trace of $t$ is negative, the rank of $\lpt$ is $2^{d-2}$ and the fixed point sublattice is a scaled version of $\bw {d-2}$, whose automorphism group is $\brw {d-2}$ if $d \ne 5$ and is $W_{E_8}$ if $d=5$.  
\end{thm}

\begin{thm}\labtt{maintheorem5}
The automorphism groups of the involution fixed point sublattices is inherited when
 the involution is dirty, split, of defect is 1 and when 
 $d\ge 5$ is odd.  
 \end{thm}
 
 \begin{thm}\labtt{maintheorem6}
The automorphism groups of the involution fixed point sublattices is not inherited when
 the involution is nonsplitsplit, of defect is 1 $d\ge 5$.  The fixed point sublattices are isometric to $ss \bw {d-2}\perp ss\bw {d-2}$.  
 \end{thm}

The author thanks Alex Ryba for many useful discussions.  
The author has been supported by
NSA grant USDOD-MDA904-03-1-0098.

\section{\bf Notation and terminology }

We mention some special terminology, definitions  and notation; see \cite{bwy}.   

\bigbreak
\halign{#\hfil&\quad#\hfil\cr

$\bw d$, the Barnes-Wall lattice in dimension $2^d$ & \cite {bwy} \cr

 $BRW^0(2^d,\pm )$ & Bolt, Room and Wall group, \cite{bwy}   \cr

clean & an element of    $BRW^0(2^d,\pm )$\cr
&  not conjugate to its negative   \cr

$D$, a lower dihedral group  & a dihedral group of order 8 \cr
& in the lower group $R$ \cr 

defect of an involution & \reftt{defect} \cr

density, commutator density & \cite{bwy}  \cr 

 determinant of a lattice, $L$    & $|\dg L |$ 
\cr 

diagonal & \reftt{diagonalgroup} \cr 

dirty & an element of    $BRW^0(2^d,\pm )$ \cr 
&  conjugate to its negative  \cr 

$\dg L$, discriminant group of an integral lattice $L$ & $\dg L = \dual L / L$\cr

$\dual L$, the dual of the lattice $L$ & $\{ x\in \QQ \otimes L | (x,L)\le \ZZ \}$ \cr

$\vep_S$ & \reftt{diagonalgroup} \cr 

fourvolution &  a linear transformation \cr 
& whose square is $-1$  \cr

$G=\gd d$ & $\brw d$ \cr 

inherited & \reftt{inh} \cr 

lower & in $R$  \cr

$R =  \rd d$ & $O_2(\brw d )$ \cr

SSD, semiselfdual,  
RSSD, relatively semiselfdual   & applies to certain sublattices \cr 
& of an integral lattice; \cr 
& there are associated involutions, \cr 

sBW, $s\bw k$   & scaled copy of some $\bw k$ \cr
& $\cong \sqrt s \bw k$ for some integer $s>0$.\cr

ssBW, $ss\bw d$  (for a sublattice of $\bw d$) & suitably scaled copy of $\bw k$ = \cr
& a scaled $\bw k$ with scale  \cr
& $2^h$, $h=\frac {d-k} 2$ for $d-k$ even;  \cr 
& $h=\frac {d-k-1}2$ for $d-k$ odd, $d$ even; \cr
& $h=\frac {d-k-1}2+1$ for $d-k$ odd, $d$ odd. \cr

total eigenlattice, $Tel(E), Tel(L,E)$ & the sum of the eigenlattices of \cr & an elementary abelian 2-group \cr
& or involution $E$ on  the lattice $L$  \cr

upper & in $G \setminus R$  \cr

\cr }

{\bf Conventions.  }  Our groups and most endomorphisms act on the
right, often with exponential notation.  
Group  theory notation is
mostly consistent with \cite{Gor, Hup, G12}. 
The commutator of $x$ and $y$ means $[x,y]=x^{-1}y^{-1}xy$ and
the  conjugate of of
$x$ by
$y$ means
$x^y:=y^{-1}xy=x[x,y]$.  These notations extend to actions  of a group on
an additive group.  

Here are some fairly standard 
notations used for particular extensions of
groups: 
$p^k$ means an elementary abelian $p$-group;
$A.B$ means a group extension with normal 
 subgroup $A$ and quotient $B$;  
$p^{a+b+\dots }$ means an iterated group extension, with factors $p^a, p^b,
\dots $ (listed in upward sense); 
$A{:}B, A{\cdot}B$ mean, respectively, a  split extension, nonsplit
extension.

\section{Preliminaries}

\subsection{Groups}

\begin{de}\labttr{evenodd} 
The {\it Dickson invariant} is a natural homomphism $O^+(2d,2) \rightarrow \ZZ_2$ which has the property that it is nontrivial on orthogonal transvections.  (For an exact definition, see \cite{Dieud}).  The kernel is  the subgroup $\Omega^+(2d,2)$.  Elements of the latter group are called {\it even} and elements of $O^+(2d,2)$ which are not even are called {\it odd}.  

This notion extends to the full holomorph
$2^{1+2d}.O^+(2d,2)$ in $GL(2^d,\CC )$, so that the BRW group $\brw d$ is considered its even subgroup  \cite{GrMont}.   
\end{de} 

\begin{nota}\labttr{gdrd}  From now on, $d\ge 2$, $\gd d:=\brw d$, $\rd d:=O_2(\gd d)$.   Reference to $d$ will typically be suppressed  and we use $G$ for $\gd d$ and $R$ for $\rd d$.  
\end{nota} 

\begin{lem}\labtt{commabel}
Let $t$ be an isometry of $V$, a vector space in characteristic 2 with an alternating bilinear form.  Then $[V,t]=Im(t-1)$ is totally isotropic.  
\end{lem}
\pf
Let $x, y \in V$.  Then $(x(t-1),y(t-1))=(x,y)-(x,yt)-(xt,y)+(xt,yt)$.  Since we are in characteristic 2 and $t$ is an isometry, the first and last terms cancel.  Since $t^2=1$, the middle two terms cancel.
\eop

\begin{rem}\labttr{commabel2} When $t$ leaves invariant a quadratic form associated to the alternating bilinear form, the totally isotropic space of \ref{commabel} may be totally singular or not. 
\end{rem}

\begin{nota}\labttr{cmz} Let $R$ be an extraspecial group and $H$ a subgroup of $R$ which contains $Z(R)$.  Then $H$ has a central product decomposition, $H=AB$, where $A=Z(H)$ and $B=Z(R)$ or $B$ is extraspecial.  Clearly, $A \cap B=Z(R)$.  The group $B$ is not unique if $A>Z(R)$, but the set of such $B$ forms an orbit under $Stab_{Aut(R)}(H)$ if $A$ is elementary abelian.  We call such a decomposition of $H$ a {\it CMZ-decomposition} (for complement modulo the center) and such a $B$ is called a {\it CMZ-subgroup}.  
\end{nota}

\begin{lem}\labtt{fixinvol}  An involution $t$ which acts on an extraspecial group $R \cong 2^{1+2d}_+$ as an even automorphism fixes a noncentral involution if $d\ge 2$.  
\end{lem}
\pf
If $t$ is inner, this is obvious.  Suppose that $t$ acts nontrivially on the Frattini factor of $R$.  Since $[R,t]$ is not contained in $Z(R)$ and is normal in $R$, $Z(R)\le [R,t]$.  Also, $[R,t]$ is abelian (by \reftt{commabel}).  Since $t$ inverts a set of generators for $[R,t]$, it inverts $[R,t]$, so centralizes $\Omega_1([R,t])$.  Also, $[R,t]$ is noncyclic since for even orthogonal transformations, the space of fixed points is even dimensional (see \reftt{defect}).  This completes the proof.  \eop

\begin{lem}\labtt{involonexsp} Let $t$ be an upper involution in the automorphism group of an extraspecial 2-group of plus type.  Then $t$ 
centralizes a maximal elementary abelian subgroup if and only if its image in the outer automorphism group is even and $[R,t]$ is elementary abelian.  
\end{lem} 
\pf  The necessity follows from the  well-known facts that $\Omega^+(2d,2)$ has two orbits on maximal totally singular subspaces and that they are fused by $O^+(2d,2)$ \cite{GrElAb}. 

We now prove sufficiency.  
We may assume that the order of the extraspecial group $R$  is $2^{1+2d}$, for $d\ge 2$ (there are no even upper involutions for $d=1$).  Let $t$ be an upper involution.  

The action of  $t$ fixes a noncentral involution $u\in R$, by \reftt{fixinvol}.  So, $t$ acts on $C_R(u)/\la u \ra \cong 2^{1+2(d-1)}_+$.  If $u \not \in [R,t]$, then $t$ acts evenly on this extraspecial group and we finish by induction.  Therefore, we are done if $t$ fixes an involution outside $[R,t]$, so suppose that none exist.  Then since $R$ has plus type, $[R,t]$ has order $2^{d+1}$.  Since $t$ inverts $[R,t]$, we are done since $[R,t]$ does not have exponent 4.  
\eop

\begin{prop}\labtt{equiv}
We are given $V=\FF^{2d}$ with quadratic form $q$ and associated bilinear form $(\cdot ,\cdot )$
so that $V=I\oplus J$ is a decompostion into 
maximal totally singular $d$-dimensional  subspaces.  
Define $Inv(V,I)$ to be the set of involutions $t$ in $G$, the orthogonal group for $q$, so that 
$t$ is trivial on $I$ and $V/I$ and $[V,t]=I$.  Then 

(0) $Inv(V,I) \neq \emptyset$ if and only if $d$ is even.

(1) Assume that $d$ is even.  Then $Inv(V,I)$ is in bijection with these two sets: 

(1.a) the set of $2d\times 2d$ matrices of the form $I_{2d}+N$, where $N$ has rank $d$ and is supported in the upper right $d \times d$ submatrix, which is alternating.  

(1.b) The set of all sequences $v_1, w_1, \dots , v_d, w_d$ with each $v_j\in J, w_j \in I$ so that 
$[v_i,t]=w_i$ for all $i$ and $(v_i,w_j)=0$ except for $\{i,j\}$ of the form $\{2k-1,2k\}$ for $k=1,\dots ,\frac d 2$ in which case  $(v_i,w_j)=1$.  
\end{prop}
\pf    For (0), use \reftt{defect}.  The proof of (1) is formal.  
\eop

\begin{de}\labttr{naturalbrw} A {\it natural BRW subgroup} of $G$ is a subgroup of the form $C_G(S)$, where $S$ is a plus type  extraspecial subgroup of $R$.  Natural BRW subgroups occur in pairs, each member being the centralizer in $G$ of the other.  
\end{de}

We need to discuss normalizers of lower elementary abelian subgroups in $G$ and centralizers of clean upper involutions.

\begin{prop}\labtt{normlowerelab} 
Let $E$ be a lower elementary abelian group of order $2^{a+b}$, where $2^a=|Z(R) \cap E|$.  Let $N:=N_G(E)$ and $C:=C_G(E)$.  Suppose that $b\ge 1$.  
Then $N$ and $C$ have the following structure.  

There are subgroups $S, T\le R$ and $P\le G$ so that 

(a) $T$ and $S$ are extraspecial of respective orders $2^{1+2(d-b)}, 2^{1+2b}$ (though $T=1$ if $b=d$), $[T,S]=1$ and $R=TS$;

(b) $EZ(R)$ is maximal elementary abelian in $S$; it follows that $TEZ(R)=C_R(E)$.

(c) the group $P:=C_N(C_R(E)/EZ(R))\cap N_N(E_0)$, 
where $E_0$ complements $Z(R)\cap E$ in $E$, satisfies  
$P\cap S=EZ(R)$ and 
$P/T \cong 2^{{b \choose 2}+b(2d-2b)}{:}GL(2b,2)$;

(d) $C_C(S)=C_G(S)$ is the natural $BRW$-subgroup containing $T$;  

(e) $C_P(T)S/S$ has the form $2^{{b \choose 2}}{:}GL(2b,2)$.  

(f) $C=O_2(P)C_G(S)$;

(g) if $a=0$, $N=CP$ and if $a=1$, $N=CSP$.  
\end{prop}

\begin{de}\labttr{mns} Given an involution $t$ in an orthogonal group over a field of characteristic 2, a {\it MNS-subspace for $t$} (minimal nonsingular) 
is a nontrivial, nonsingular subspace which is $t$-invariant, and no proper subspace of it has these properties.  
\end{de}

\begin{lem}\labtt{mnslist} Let $t$ be an involution in the orthogonal group  $\Omega^\vep(2e,2)$ and $S$ a MNS-subspace for $t$.  Suppose that $t$ acts nontrivially on $S$.  

Either $S$ has dimension 2 and a basis $u, v$ so that $u^t=v$ and $(u,v)=1$, so that $u$ and $v$ are both singular or both nonsingular;

 or $S$ has dimension 4 and a basis $u_1, u_2, v_1, v_2$ of singular vectors so that $v_1^t=v_2, u_1^t=u_2$ and the Gram matrix for this basis is $\begin{pmatrix} 
0&0&1&0\cr
0&0&0&1\cr
1&0&0&0\cr
0&1&0&0\cr
\end{pmatrix}$.   Furthermore both spaces are MNS-subspaces.  
\end{lem}
\pf
We may suppose that 
$dim(S)\ge 4$ and that 
for every singular vector $v \in S$, $(v,v^t)=0$, then try to get the last conclusion.  We note  that $S$ is spanned by its singular vectors.  

Take a singular vector $v_1$ not fixed by $t$ and define $v_2:=v_1^t$.  
Choose a singular vector $u_1\in S$ so that 
$(v_1,u_1)=1$ and $(v_2, u_1)=0$.  Using $t$-invariance, we find that the sequence $v_1, v_2, u_1, u_2:=u_1^t$ has Gram matrix 
$\begin{pmatrix} 
0&0&1&0\cr
0&0&0&1\cr
1&0&0&b\cr
0&1&b&0\cr
\end{pmatrix}$.  This matrix is nonsingular, whence $S$ has dimension just 4.  Now, if $b\ne 0$, $span\{v_1+u_2, u_1+v_2\}$ is a 2-dimensional MNS-subspace.  Therefore, $b=0$.  Since $S(t-1)$ is totally singular, $S$ is minimal.  
\eop

\begin{lem}\labtt{swap} 
Let $u$ be an involution in $\Omega^+(2e,2)$ of defect $e$.    There exists a maximal totally singular subspace $F$ so that $F \cap F^u=0$.  
\end{lem}
\pf 
Take a MNS-subspace for $S$.  Then $t$ acts nontrivially on $S$ since the defect is $e$.  
Also, $t$ leaves invariant the summands of the decomposition $S \perp S^\perp$.  We are therefore done by induction if we check it for the cases of \ref{mns}.  This is trivial for the 2-dimensional case and for the 4-dimensional case, take the span of the second and third basis elements.  
\eop

\begin{nota}\labttr{diagonalgroup}  
On the rational vector space spanned by a Barnes-Wall lattice, we take a sultry frame $F$ containing a basis labeled by affine space $\FF_2^d$ \cite{bwy}.  For a subset $S$ of the index set, define the orthogonal  involution $\vep_S$ to be the map which is $-1$ at frame elements labeled by a member of $S$ and 1 on the other frame elements.  The set of such linear maps, for $S\in RM(2,d)$, forms the {\it diagonal group}, denoted $\cal E$ or ${\cal E}_d$.  It is a subgroup of $\brw d$.  The defect of the codeword $c$ is the defect of the involution $\vep_c$.   
\end{nota}

\begin{de}\labttr{cleandirty2} 
Recall that an involution in the BRW group $\brw d$ 
is {\it dirty} if 
it is conjugate to its negative and otherwise, it is {\it clean}; \ref{cleandirty}.  These properties are equivalent to having nonzero, zero trace, respectively, on the natural $2^d$-dimensional module.  Furthermore, if the trace is nonzero, it has the form $\pm 2^{d-k}$, where $k$ is the defect \reftt{defect} of the involution.  
We call such an involution a {\it $(d,k)$-involution}.  Any involution in the lower coset of such is also called a {\it $(d,k)$-involution}.  

The dimension of the space of commutators of a defect $k$ diagonal involution with the translation group of $AGL(d,2)$ is $2k$ since the translation group can be interpreted as a complement in $\rd d$ to the diagonal subgroup corresponding to $RM(1,d)$.    
The terms clean and dirty apply to codewords, according to whether the corresponding involutions are clean or dirty.

The term {\it absolute clean trace} or 
{\it positive clean trace} applies to any element of $\brw d$ and means, the absolute value of the trace of any clean element in its lower coset.  So, the absolute clean trace is a power of 2 even if the element is dirty.  
We let $\cal D$ and $\cal C$, respectively, denote the set of dirty and clean codewords in $RM(2,d)$.  
\end{de}

\begin{prop}\labtt{centupperinvol} 
Let $u\in G$ be a clean $(d,k)$-involution, $k>0$.  Then 

(i) $C_G(u)$ has the following form: it is a subgroup of $N_G(E)$, where $E=[R,u]$ is a rank $2k+1$ elementary abelian group as in \reftt{normlowerelab}; $C_G(u)$ corresponds to the natural $Sp(2k,2)$ subgroup of $GL(2k,2)$  associated to the identification of $R/C_R(E)$ with $E/Z(R)$ derived from commutation with $t$;

(ii) The involution $uR\in G/R$ has centralizer $C_G(u)R/R$.  
\end{prop}  
\pf 
(i) 
It is clear from \reftt{normlowerelab} that $C_G(u)$ has this form,   except possibly for the replacement of $GL(2k,2)$ by $Sp(2k,2)$.  It is clear that commutation  by $u$ gives a linear isomorphism of $S/E$ onto $E/Z(R)$ which makes these two spaces into dual modules for $C_G(u)$.  The action of $C_G(u)$ is therefore symplectic on both.  It suffices to show that 
there is a subgroup of $C_G(u)$ which acts on both as the full group $Sp(2k,2)$.  

We take an elementary abelian subgroup $F$ of $S$ so that 
$FZ(R)=F\times Z(R)$ is maximal elementary abelian and so that 
$F \cap F^u=1$ (see \reftt{swap}).   Then $u$ acts on $H:=C_{C_G(u)}(T)\cap N_G(F) \cap N_G(F^u)$, which has shape $2\times GL(2(d-k),2)$ (the shape is clearly of the form $2.GL(2k,2)$ but is actually a direct product; see \cite{bwy} or the Appendix).   
Clearly, $C_H(u)$ has shape $2 \times Sp(2k,2)$.  

(ii) This follows from noticing that the set of clean elements in $uR$ is just the union of the $R$-conjugacy class of $u$ with the $R$-conjugacy class of $-u$.  
\eop 

\begin{rem}\labttr{centdirty} The exact 
 structure of centralizers for dirty involutions is not needed in this article, but we give a sketch.  

There are three main kinds of dirty involutions: lower involutions (defect 0); upper split (positive defect, with elementary abelian commutator subgroup on $R$); (upper)  nonsplit (positive defect, with exponent 4 commutator subgroup on $R$).   

The centralizer of a lower involution has shape $[2\times 2^{1+2(d-1)}]2^{2(d-1)}.\Omega^+(2(d-1),2)$.  

Let $t$ be a dirty split upper involution.  Then $t=ru$, where $u$ is an upper involution 
and $r$ is a lower involution from $R\setminus [R,u]$.  The structure of $C_G(u)$ is discussed in \ref{centupperinvol}.  We have $C_G(t) \le C_G(u)$, $C_R(t)$ has index 2 in $C_R(u)$ and $C_G(t)R/R$ is a natural subgroup of $C_G(u)R/R$ of shape $2^{2(d-2k)}{:}\Omega^+(2(d-2k),2)$.  

Let $t$ be a nonsplit involution.  Let $S$ be a maximal extraspecial subgroup of $C_R(t)$.  Then $C_R(S)\ge [R,t]=[C_R(S),t]$.  
Also, $C_R(t)=S \times E$, where $E$ is elementary abelian and a complement in $\Omega_1([R,t])$ to $Z(R)$.  
We say $t$ has {\it plus type} or {\it minus type} according to the type of the extraspecial group $S$.  Now, $N_G([R,t])\ge R$ and 
$N_G([R,t])/R$ modulo its unipotent radical has the form $\Omega^+(2(d-2k),2) \times GL(2k-1,2)$.  
The image of $C_G(t)$ in the latter quotient has the form 
$\Omega^+(2(d-2k),2) \times O(2k-1,2)$. 
\end{rem}

\subsection{The codes $RM(2,d)$ and the diagonal group}

Our vector spaces are finite dimensional.  We shall mix styles at times, so that a codeword may be written in lower case (when we think of it as a vector) or upper case (if we think of it as a geometric structure, like an affine subspace).

\begin{nota}\labttr{rm}  The {\it Reed-Muller code} $RM(k,d)$ is the binary code indexed by affine space $\FF_2^d$ and spanned by all affine subspaces of codimension $k$.  Its dimension is 
$\sum_{i=0}^k {d \choose i}$. 
\end{nota}

\begin{de}\labttr{midset} 
A {\it midset}  is a codeword in $RM(2,d)$ of size $2^{d-1}$.  A midset is {\it nonaffine} if it is not a codimension 1 affine subspace.  
A codeword is {\it short} if its weight is less than $2^{d-1}$.  
A codeword is {\it long} or {\it tall} if its weight is more than $2^{d-1}$.  
\end{de}

\begin{lem}\labtt{involutionisdiagonal}
Let $t \in G$ be an involution so that $[R,t]$ is elementary abelian and $\cal E$ a given diagonal group. Then there is a conjugate of $t$ in $\cal E$,  unless possibly $d$ is even and $t$ has defect $\frac d 2$, in which case there exists another diagonal group containing $t$.
\end{lem}
\pf
Use \ref{involonexsp} and the fact that $C_R(t)$ is nonabelian if and only if $C_R(t)$ contains representatives of both $G$-conjugacy classes of maximal elementary abelian subgroups of $R$.   
\eop

\begin{nota}\labttr{aglnota} 
We will study the action of $AGL(d,2)$ on $\FF_2^d$ and various codes.  Let $T:=T(d,2)$ denote the translation subgroup and $GL(d,2)$ the stabilizer of some origin (understood from context).  
\end{nota}

\begin{de}\labttr{coindep}   
Linear subspaces $U_i$ of a vector space are {\it independent} if their sum is their direct sum.  
Linear subspaces $U_i$ of a vector space are {\it coindependent} if their annihilators in the dual space are independent. 
 
This definition extends to a collection of affine subspaces $U_i$ of a vector space, provided their common intersection is nonempty.  One then chooses any origin in $\bigcap_i  U_i$ and uses the above definition (which is independent of choice of origin).  
\end{de}

\begin{lem}\labtt{kfoldint2} Suppose that we have $k\ge 1$ linearly coindependent codimension 2 affine subspaces $S_1, \dots , S_k$ in $\FF_2^d$ with nonempty common intersection.  Then $|S_1+ \dots +S_k|= 2^{d-1}-2^{d-k-1}$.  (Note: $k\le \frac d 2$ here.) 
\end{lem}
\pf
Let $a(d,k)$ be $ 2^{d-1}-2^{d-k-1}$.  
We use induction on $k$.   The result is trivial for $k=1, 2$.  We may assume that the spaces contain a common origin, so are linear.  

  Assume that $k\ge 3$ and that the formula holds by induction for $k-1$.  We have $S_k \cap (S_1+\dots + S_{k-1})=S_1\cap S_k + \dots + S_{k-1}\cap S_k$, which, by induction on $d$ and coindependence in $S_k \cong \FF_2^{d-2}$, has cardinality $a(d-2,k-1)$.  It follows that $|S_1+ \dots +S_k|= 2^{d-2}+ a(d,k-1)-2a(d-2,k-1)=a(d,k)$.   
\eop  

\begin{de}\labttr{cubi} A set of codimension 2 subspaces as in \reftt{kfoldint2} is called a {\it cubi sequence of codimension 2 spaces}.   
Their Boolean sum is called a {\it a cubi sum}.  
\footnote{We chose the term {\it cubi} because our theory suggested  the remarkable cubi sculpture series by David Smith.  See also the footnote at \ref{cleansingexample}.  
} 
\end{de}

\begin{nota}\labttr{orderedcubis}  Let $c$ be a clean codeword of defect $k$.  Let $$Cubi(c):=\{(S_1,\dots ,S_k)|\bigcap_{i=1}^k S_i \ne \emptyset,  S_1,\dots ,S_k \hbox{ are coindependent affine }$$ 
$$\hbox{ codimension 2 subspaces, and }\sum_{i=1}^k S_i=c \},$$
the set of {\it cubi expressions} of $c$, i.e. the set of ordered cubi sequences as above whose sum is $c$. 
\end{nota}

\begin{coro}\labtt{alldefects}
Given any integer $j \in [0, \frac d 2]$, there is an involution of defect $j$ in the diagonal group. 
\end{coro}
\pf  
If $j=0$, take a lower involution.  Suppose $j>0$.  
Then take $\vep_{S_1+ \dots +S_j}$, in the notation of \reftt{kfoldint2}.  
\eop

Next, we show explicitly how to realize a dirty class associated to the clean class within the diagonal group.  

\begin{lem}\labtt{cleantodirty} 
Given $d\ge 3$ and $k\ge 1$ and a length $k$ cubi sequence 
in $\FF_2^d$, there exist hyperplanes whose sum with the cubi sum has cardinality $2^{d-1}$.  
In fact, any hyperplane 
which neither  contains nor avoids the cubi intersection meets this condition.  
\end{lem}
\pf
Let $S_1, \dots , S_k$ be our cubi sequence and let $U:=\bigcap_{i=1}^k S_i$.  Let $\cal N$ be the set of hyperplanes  which neither contain $U$ nor avoid $U$.   Then $|{\cal N}|=2^{d+1}-2^{2k+1}$.   This is positive for $d\ge 3$ and $k \ge 1$.  

Let $H \in {\cal N}$.  
Then the spaces $S_i \cap H$ have codimension 2 in $H$.  They are coindependent with respect to $H$ since $H\cap U$ has codimension 1 in $U$.  
Therefore, \reftt{kfoldint2} gives 
$|H \cap (S_1+\dots +S_k)|=|(S_1\cap H) +\dots +(S_k\cap H)|=2^{d-2}-2^{d-k-2}$.  Consequently, 
$|H+S_1+\dots +S_k|=2^{d-1}+2^{d-1}-2^{d-k-1}-2(2^{d-2}-2^{d-k-2})=2^{d-1}$.  
\eop

\begin{rem}\labttr{nothyp}  The codeword of weight $2^{d-1}$ constructed in the proof of \reftt{cleantodirty} is not a hyperplane, since the Boolean sum of two distinct nondisjoint hyperplanes is a hyperplane and $|S_1+\dots +S_k|<2^{d-1}$.  
\end{rem}

We next need to work from a nonaffine midset to the class of clean codewords that it comes from.

\begin{de}\labttr{cleansing}  Let $d\ge 3$.  
Given a nonaffine midset $a$, a hyperplane $h$ so that $a+h$ is clean is called a {\it cleansing hyperplane for $h$}.  It follows that if $a$ has defect $k$, and $h$ is cleansing, then $|a\cap h|=2^{d-2}\mp 2^{d-k-2}$.  (Note that $d-k\ge 2$ for $d\ge 3$.) 
\end{de}

\begin{lem}\labtt{cleanrep} Every coset of $RM(1,d)$ in $RM(2,d)$ contains a clean codeword.  
\end{lem}
\pf
Take a nontrivial coset, say $u+RM(1,d)$ and take a complement $S$ in $RM(1,d)$ to the 1-space spanned by the universe.  
The subgroup of the diagonal group corresponding to $S$ has 1-dimensional fixed point sublattice, so the sum of the traces of its elements is $2^d$.  Assume that the lemma is false.  Then every element of  $\la u, S \ra \setminus S$ gives a diagonal map of trace 0. Therefore the sum of the traces for 
the subgroup of the diagonal group corresponding to $\la u, S \ra$ 
is $2^d$, which is impossible since this number must be divisible by $2^{1+d}$.  
\eop

\begin{lem}\labtt{transforms} 
If $c \in RM(2,d)$ is clean, the number of its conjugates by $R$ is $2^{2k}$, where $c$ has defect $k$. 
\end{lem}
\pf 
This is just the correspondence of the $R$-orbit of $c$ under the action of conjugation on $RM(2,d)$ with the cosets of $C_R(c)$ in $R$, together with the definitions of defect and cleanliness.  
\eop

\begin{prop}\labtt{cosetcleanliness} In a given coset $c+RM(1,d)$, where $c$ is clean and has defect $k$, the number of clean codewords is $2^{2k+1}$ and the number of dirty codewords is $2^{d+1}-2^{2k+1}$.  
\end{prop}
\pf  
If $c \in RM(2,d)$, the number of its transforms by $R$ is $2^{2k}$, by  \reftt{transforms}.  
The coset $c+RM(1,d)$ also contains the same number of transforms of the complement $c+\FF_2^d$, which is also clean.  

We use the irreducible module for $G$,  which is a 
 $2^d$-dimensional complex vector space, and the trace function $Tr$ on it.  
The previous paragraph implies that the sum 
$s(c):=\sum_{v \in c+RM(1,d)} Tr(v)^2$ is at least $2\cdot 2^{2k+2(d-k)}=2^{2d+1}$

Since the group $RM(1,d)$ acts on the $2^d$-dimensional complex vector space so as to afford all linear characters nontrivial on the center, each with multiplicity 1, it follows from orthogonality relations for the group generated by $R$ and $c$ that each $s(c)=2^{2d+1}$.  The coset therefore has $2^{2k+1}$ clean elements and $2^{d+1}-2^{2k+1}$ dirty elements.  
\eop

\begin{coro}\labtt{numcleansing}
The number of cleansing hyperplanes for a dirty codeword $s\in RM(2,d)$ is $2^{2k+1}$, where $k$ is the defect of any clean involution in the coset  $s+RM(1,d)$.  Thus the set $\cal N$ of \reftt{cleantodirty} is the full set of noncleansing hyperplanes.   
\end{coro}

\begin{ex}\labttr{cleansingexample}  Let $d=4, k=1$ and let $S$ be a defect 1 (nonaffine) midset.   There are 8 cleansing hyperplanes.  Write $S=A+H$, where $A$ is short and $H$ a cleansing hyperplane of $S$ (this involves half the cleansing hyperplanes).  Then $A$ is a 4-set (hence an affine hyperplane) and $S\cap H$ is a 2-set.  This set is stable by translation with elements of the core.  Therefore, $S$ is a union of four cosets of $S\cap H$.  The assignment $H\mapsto S\cap   H$ is one-to-one from the set of cleansing hyperplanes such that $S+H$ is short. 
By counting, this is a bijection.  The union of any two sets $S\cap H$, as $H$ varies, is an affine 2-space.  Therefore, $S$ is the disjoint union of a pair of disjoint, nonparallel affine 2-spaces, in three different ways.  
\footnote{These configurations also suggest the David Smith cubi theme; see \ref{cubi}.  }
\end{ex} 

\begin{coro}\labtt{intersectcleansing} Given cleansing hyperplanes $H_1, H_2$ for the dirty codeword $S$, if $H_1\cap S=H_2\cap S$, then $H_1=H_2$, i.e.,  for cleansing hyperplanes, $H$, the map $H\mapsto H \cap S$ is monic.  
\end{coro}
\pf
If $H_1$ and $H_2$ are distinct, then, since they meet, their sum is a hyperplane.  Since $H_1+H_2 $ is contained in the complement of $S$, it equals the complement of $S$.  This is a contradiction since $S$ is not affine.  
\eop

\begin{proc}\labttr{proc} We now have a {\it procedure to determine the orbit of a dirty codeword}.  It depends only on examining the code, not the action of the group $AGL(d,2)$.  
Call such a codeword $v$.  Add to $v$ all of the $2^{d+1}-2$ affine hyperplanes.  A nonempty set of these will be cleansing and the corresponding sums will have weight of the form $2^{d-1}\pm 2^{d-k-1}$, which will give 
the defect $k$.  This procedure is exponential in $d$.  
\end{proc}

\begin{lem}\labtt{localcleanl}
Two short (resp. long) clean codewords of the same defect are in the same orbit under $AGL(d,2)$.  A short clean codeword is a cubi sum.  
\end{lem}
\pf  We interpret these codewords by their actions on the commutator quotient of $R$.  The result follows from transitivity of the natural action of $GL(d,2)$ on alternating matrices of the same rank.
\eop

\begin{lem}\labtt{indepofcubi} Suppose that we are given $(S_1,\dots ,S_k)\in Cubi(c)$ as in \reftt{orderedcubis}.    The subspace $\bigcap_{i=1}^k S_i$ has dimension $d-2k$ and is the subgroup of the group of translations which fixes $c$.   This subspace depends on $c$ only, not on a choice from $Cubi(c)$.  
\end{lem}
\pf
Clearly, the above intersection is a linear subspace and translations by it fix each $S_i$, hence also fix $c$.  Since the space of commutators of the translation group with $c$ has dimension $2k$, no translations outside this subspace fixes $c$.  Therefore, this intersection depends on $c$ only.  
\eop

\begin{lem}\labtt{transobclean} The stabilizer in $AGL(d,2)$ of the clean codeword $c$ of defect $k$ is transitive on $Cubi(c)$, and the the stabilizer of a member of $Cubi(c)$ has shape $2^{d-2k}.2^{2k(d-2k)} [(\prod_{i=1}^kGL(2,2)) \times GL(d-2k,2)]$.  
\end{lem}
\pf  The initial $2^{d-2k}$ refers to the group of translations which stabilize $\bigcap_{i=1}^k S_i$.  The result follows from transitivity of $GL(d,2)$ on ordered direct sums of $k$ 2-spaces in the dual.  
\eop

\begin{de}\labttr{core} The {\it core} of a clean codeword is  $\bigcap_{i=1}^k S_i$, where $(S_1,\dots  ,S_k)\in Cubi(c)$.  
The definition is independent of choice from $Cubi(c)$, by \reftt{indepofcubi}.  
\end{de}

\begin{thm}\labtt{cleancount} 
The stabilizer of a clean codeword of defect $k$ in $AGL(d,2)$  is a group of the form $[2^{(1+2k)(d-2k)}]{:}[Sp(2k,2) \times GL(d-2k,2)]$.  
It has two orbits on $\FF_2^d$, namely the core and its complement.    
\end{thm}
\pf
The second statement follows from the structure of the stabilizer, which we now discuss.   
 
We may think of our clean codeword $c$ as a cubi sum for cubi sequence $(S_1, \dots ,S_k)$.  Choose an origin in  the core \reftt{core}, i.e., the $(d-2k)$-space $U:= S_1 \cap \dots \cap S_k$. 

Let $H$ be the stabilizer of $c$ in $AGL(d,2)$.   
Then $H_t:=H\cap T$ is transitive on $U$.  
The last paragraph implies that  $H=H_t H_0$ where  $H_0$ is the stabilizer of the origin.  So, $H_t$ corresponds to $U$ and $H_0$ lies in the stabilizer in $GL(d,2)$ of the subspace $U$, a parabolic subgroup $P$ of the form $2^{2k(d-2k)}{:}[GL(2k,2)\times GL(d-2k,2)]$.  Note that $O_2(P)$ is a tensor product of irreducibles for the two factors, so is irreducible.  

We next argue that $H_0$ is a natural $2^{2k(d-2k)}{:}[Sp(2k,2)\times GL(d-2k,2)]$-subgroup of $P$.

Consider $C_G(t)$, where $t$ is the diagonal matrix $\vep _c$.  Then we have the CMZ decomposition \reftt{cmz} for $C_R(t)$ and a related one for $R$:  $R=R_1R_0$, where $[R_0, R_1]=1$,  $C_R(t)=C_1R_0$, where $R_0$ is extraspecial, and $C_1 \le R_1$ and $C_1$ is elementary abelian and contains $Z(R)$.  
There is a corresponding product $J_0J_1$ of commuting natural BRW subgroups, with $R_i=O_2(J_i), i=1,2$.  
We have $|C_1|=2^{2k+1}$ and $C_1=[R,t]=[R_1,t]$.  The action of $t$ preserves $R_1$ and the maximal elementary abelian subgroup $C_1$.  Also, $t$ acts on $N_{J_1}(C_1) \cong 2^{1+4k}2^{2k \choose 2}GL(2k,2)$.  There is a pair of maximal elementary abelian subgroups $B_1, B_2$ so that $R_1=B_1B_2, B_1 \cap B_2=Z(R)$ and $t$ interchanges $B_1$ and $B_2$ (see \reftt{swap}).

Choose $D_i \le B_i$ so that $B_i=D_1 \times Z(R)$ and $t$ interchanges $D_1$ and $D_1$.  The common stabilizer of $D_1$ and $D_2$ in $Aut(R_1)$ has the form $2 \times GL(2k,2)$.  The action of $t$ has fixed point subgroup of the form $2 \times Sp(2k,2)$ because $D_1$ and $D_2$ are in $t$-invariant duality.  Therefore, the image of $H$ in the left factor of $P/O_2(P) \cong GL(2k,2)\times GL(d-2k,2)$ contains a copy of $Sp(2k,2)$.  Since the image of $H$  in the left factor stabilizes a nondegenerate form, the image is exactly $Sp(2k,2)$.  

We claim that the stabilizer of $c$ in $AGL(d,2)$ contains the natural $GL(d-2k,2)$ subgroup which commutes with the above copy of $Sp(2k,2)$.  This follows since the stabilizer of a member of $Cubi(c)$ involves a copy of $GL(d-2k,2)$ which acts faithfully on the core and commutes with the action of the above $Sp(2k,2)$, which acts trivially on the core and faithfully on a complement to the core (meaning, on a linear complement, assuming the origin is chosen from the core).

The claim implies that  $H$ maps onto 
the right factor of 
 $P/O_2(P) \cong GL(2k,2)\times GL(d-2k,2)$.  
It follows that $O_2(P)$ is an irreducible module for $H$ (a tensor product of irreducibles for the factors $Sp(2k,2)$ and $GL(d-2k,2)$), whence $H \cap O_2(P)$ is either 1 or $O_2(P)$.  The latter group preserves all cosets of $U$ in $\FF_2^d$ and each $S_i$ is a union of such cosets, whence $O_2(P)\le H$. 
\eop

\begin{lem}\labtt{localdirty}
Two dirty codewords of the same defect are in the same orbit under $AGL(d,2)$.  
\end{lem}
\pf This is obvious from \reftt{cleantodirty} and how the stabilizer of the core in $AGL(d,2)$ acts on $\FF_2^d$.  \eop 

\begin{rem}\labttr{summary} The main theorems \reftt{maintheorem1} and \reftt{maintheorem2} follow from \reftt{localcleanl} \reftt{localdirty} \reftt{liftclass}, \reftt{nonsplitinvols}.  
Note that we get as a corollary the well-known result that the minimum weight codewords in $RM(2,d)$ are the affine codimension 2 subspaces.  
\end{rem}

\begin{prop}\labtt{dirtystab} 
Let $c$ be a clean codeword of defect $k$.  

(i) The stabilizer in $AGL(d,2)$ of the coset $c+RM(1,d)$ is $T(d,2)S$, where $T(d,2)$ is the full translation group and $S$ is the stabilizer of $c$ in $AGL(d,2)$ (see \reftt{cleancount}).  

(ii) Let $s\in c+RM(1,d)$ be a dirty codeword.  
The commutator space $[T(d,2),s]$ has dimension $2k$
The stabilizer of $s$ in $AGL(d,2)$ is a subgroup of $S$ of index $2^{d+1}-2^{2k+1}$ of shape $[2^{(1+2k)(d-2k-1)}][Sp(2k,2) \times AGL(d-2k-1,2)]$.  It is $Stab_S(h)$, where $h=s+c$ is an affine codimension 1 subspace which meets the core of $c$ in a codimension 1 subspace of it.  The initial $2^{1\cdot (d-2k-1)}$ corresponds to translations by the intersection of the core of $c$ with a cleansing hyperplane.  
\end{prop}
\pf
(i) This is clear since the set of clean elements in $c+RM(1,d)$ is just the set of $2^{2k}$ $T(d,2)$-transforms of $c$.  

(ii) Since $s$ is dirty, $d-2k>0$.  

Consider the set $\cal P$ of all pairs $(s,r) \in c+RM(1,d)$ so that $s$ is dirty, $r$ is short and clean (whence $s+r$ is a hyperplane, so is a cleansing hyperplane; \reftt{cleansing}).  
We refer to \reftt{cleancount}.  
Let $H$ be the stabilizer of this coset in $AGL(d,2)$.  Then $H$ acts transitively on $\cal P$, which has cardinality $(2^{d+1}-2^{2k+1})2^{2k}$, so $Stab_H((s,r))$ has index $2^{d+1}-2^{2k+1}$ in $Stab_H(r)$, which has form 
$[2^{(1+2k)(d-2k)}]{:}[Sp(2k,2) \times GL(d-2k,2)]$.

Now, consider a hyperplane $h$
 in $\FF_2^d$ which meets $U$ in a codimension 1 subspace of $U$.  By \reftt{cleantodirty}, $r+h$ is a midset, so $(r+h,r) \in {\cal P}$.  Since $H_{(r+h,h)}$ stabilizes $h$, 
it follows that  $H_{(r+h,r)}$, hence every $H_{(s,r)}$, has the form 
$[2^{(d-2k-1)+2k(d-2k)}]{:}[Sp(2k,2) \times AGL(d-2k-1,2)]$.  
\eop

\section{The conjugacy classes of involutions in $\gd d$ and orbits on RSSD sublattices}

We continue to  let $G:=\gd d$, $R:=\rd d$ and let $t \in G$ be an involution.  We summarize the conjugacy classes of involutions.

Suppose that  $t$ centralizes a maximal elementary abelian subgroup (so is in a diagonal group).  
For each maximal elementary abelian subgroup $E$ of $C_R(t)$, we have 
representatives of 
$\lfloor \frac d 2 \rfloor$ clean classes of upper involutions in a diagonal group $C_G(E)$.  
Upper involutions of the same defect and trace are conjugate in $G$ except for the case where $d$ is even and the involutions have full defect $\frac d 2$.  Two such involutions are clean and are conjugate if and only if their traces are equal and maximal elementary abelian subgroups in their lower centralizers are in the same orbit under the even orthogonal group.

Suppose that $t$ does not centralize a maximal elementary abelian subgroup.  Then $[R,t]$ is abelian of exponent 4 and has order $2^{1+2k}$ for some $k\ge 1$.  
It is now clear 
from \reftt{liftclass} \reftt{nonsplitinvols} 
that $t$ is conjugate to some $\eta_{2k,\pm}$ \reftt{eta}.

\begin{proc}\labttr{distinguishdirty}  
In \cite{bwy}, we showed that two RSSD sublattices in $\bw d$ which had the same rank, but unequal to $2^{d-1}$ (the clean case), are in the same orbit under $\brw d$ with the exception of two orbits for maximal defect $\frac d 2$.  
Also, \cite{bwy} treats the case of rank $2^{d-1}$ sublattices which are fixed points of lower involutions.   We now give a procedure for determining when two RSSD sublattices are in the same orbit of $\brw d$ {\it which depends only on examining a restricted set of sublattices, not the whole group $\brw d$}.  Besides the two given RSSD sublattices, we need to examine only the ones associated to lower involutions, which may be constructed directly, by induction.  

Recall that for $d>3$,  the lower involutions in $\brw 3$ are those RSSD involutions associated to $ss\bw {d-1}$ sublattices \cite{bwy}.  

Here we deal with the general dirty case, i.e.,   rank $2^{d-1}$, which represents many orbits.    Their associated RSSD involutions are dirty, so if diagonalizable are conjugate to elements of the diagonal group supported by a midsize codeword.   We assume that $d>3$.

We are given a dirty RSSD sublattice.  
Multiply this involution by all lower 
involutions.

Suppose that 
a nonempty set of such products are clean involutions with common defect $k \in [0,{\frac d 2}]$.  Since the defect $k$ is less than $\frac d 2$, $k$ determines the orbit of the sublattice, by \reftt{dirtystab}.   If $k=\frac d 2$, there are two orbits, depending  on which maximal elementary abelian lower group corresponds to the RSDD involution.  

Suppose that no such product is clean.  Then the involution is some $\eta_{2k,\pm}$.  
The subgroups $C_R(t)$ and $[R,t]$ determine $k$ and the sign $\pm$ and so the orbit of the sublattice.  
\end{proc}

For completeness, we treat the case $d=3$.

\begin{prop}\labtt{orbitsd=3} In $\bw 3 \cong \leh$, the orbits of $\weh$ on RSSD sublattices are (i) those of $\brw 3$ on RSSD sublattices of even rank, i.e., one for rank 2, two for rank 4 and one for rank 6; and (ii) four orbits, of respective ranks 1, 3, 5, 7, which are sublattices generated by a root, a set of three orthogonal roots, and the annihilators of such sublattices.
\end{prop}
\pf 
Note that the determinant 1 subgroup of $\weh$ contains a natural $\brw 3$ subgroup of odd index.  For  rank 2 and 6 sublattices, we are in the clean cases in $\brw 3$.  For rank 4, we are in the dirty cases, of which there are just two, associated to  a nonsplit involution and to a lower involution.   (There are no upper dirty involutions for $d=3$.)

There are two orbits of $\weh$ on 4-sets of mutually orthogonal pairs consisting of roots and their negatives.  One of these 4-sets spans a sublattice of $\bw 3$ which is a direct summand and the other spans a sublattice contained in a $D_4$-sublattice.  These cases correpond in the above sense to the nonsplit and lower cases.  

Now consider the case of odd rank fixed point sublattice, $M$.  It suffices to do the ranks 1 and 3 cases.  We use a lemma that if $g$ is in a Weyl group and $V$ is the natural module, then $g$ is a product of reflections for roots which lie in $[V,r]$ \cite{Car}.  At once, this implies that the rank 1 lattice here is spanned by a root.  
Suppose now that $rank(M)=3$.  
Let $\Phi$ be the set of roots in $M$.  If there is a pair of nonorthogonal linearly independent roots, then $\Phi$ has type $A_3$ or $A_2A_1$.  Since $\dg M$ is an elementary abelian 2-group, neither of these is possible.  We conclude that $\Phi$ has type $A_1A_1 A_1$.  Since $M$ is even, it must equal the sublattice spanned by $\Phi$.  We are done since $\weh$ has a single orbit on subsets of three orthogonal roots in a root system of type $E_8$.  
\eop

\begin{rem}\labttr{smallrank} 
For simplicity, discuss the main theorems for ranks at most 3 so that we may later use the assumption $d\ge 4$, as needed.

When $d=1$, the fixed point sublattice of any involution is 0 or a rank 1 lattice.

Assume $d=2$.  The dirty involutions in $\bw 2$ and their fixed point sublattices are analyzed in \reftt{alldefectsoccur}.  If $t \in \brw 2$ is clean, its fixed point sublattice has rank 1 or 3.  In these respective cases, the sublattice is spanned by a vector of norm 2 or 4 or is the orthogonal of such a rank 1 sublattice, so is a root lattice of type $B_3$ or $C_3$.  See the proof of \reftt{orbitsd=3}.  

When $d=3$, all fixed point sublattices are accounted for in the proof of  \reftt{orbitsd=3}.  They are all orthogonal direct sums of indecomposable root lattices.  
\end{rem}

\subsection{Containments in $RM(2,d)$}

\begin{lem}\labtt{aboutcontainment}  Let $A, B \in RM(2,d)$ and suppose that $0\ne A < B \ne \FF_2^d$.  Let $X^c$ denote the complement of the subset $X$ of $\FF_2^d$.  
Then one of the following holds:

(i) $A$ is a codimension 2 subspace and $B$ is a midset; or 
$B^c$ is a codimension 2 subspace and $A^c$ is a midset.  

Furthermore, (i)  happens for affine hyperplanes $B$ for any $d\ge 3$, and for nonaffine midsets $B$ exactly when $B$ has defect $1$ and $d\ge 3$, respectively.  

(ii) $A$ is short and $B$ is long, of respective cardinalities $2^{d-1}-2^{d-k-1}$. $2^{d-1}+2^{d-r-1}$, where $(k,r)=(1,1), (1,2), (2,1)$ or $(2,2)$.  We summarize:  
$$\begin{matrix}
(k,r) & |A| & |B| & |A+B| \cr 
(1,1) & 2^{d-1}-2^{d-2}=2^{d-2} &2^{d-1}+2^{d-2}=2^{d-2}3&2^{d-1}\cr
(2,1)& 2^{d-1}-2^{d-3}=2^{d-3}3&2^{d-1}+2^{d-2}=2^{d-2}3&2^{d-3}3\cr
(1,2)& 2^{d-1}-2^{d-2}=2^{d-2}&2^{d-1}+2^{d-3}=2^{d-3}5&2^{d-3}3\cr
(2,2)& 2^{d-1}-2^{d-3}=2^{d-3}3&2^{d-1}+2^{d-3}=2^{d-3}5&2^{d-2}\cr
\end{matrix}$$

Note that cases (1,2) and (2,2) are dual in the sense that $A$ and $A+B$ may be interchanged.   
Note that the case (1,1) corresponds to (i) for the midset $A+B$ containing $B^c$.  Note also that $A$ in case (1,2) and $A+B$ in case (2,2) are codimension 2 affine spaces.

\end{lem}  
\pf 
If $B$ is a midset, and $A$ is not a codimension 2 affine subspace, then $A<B$ implies that $A$ has cardinality $2^{d-1}-2^{d-k-1}$ for an integer $k$ and $A$ is a cubi sum in the sense of  \ref{localcleanl}.  
Since $A+B=A\setminus B$ is also a codeword, it has cardinality $2^{d-1}-2^{d-r-1}$ for an integer $r \ge 1$.  It follows that $k=r=1$.  Then $A$ and $A^c$ are affine codimension 2 subspaces.  Therefore, if $B$ is a midset, (i) holds.  

It is obvious that (i) happens in an essentially unique way when $B$ is an affine hyperplane.  Assume $B$ is a midset but not affine.  The codimension 2 affine subspaces  $A$ and $A'$ whose union is $B$ are not translates of each other.  Let $A''$ be a translate of $A'$ which meets $A$ nontrivially.  
The intersection has codimension 1 or 2 in 
each of $A$ or $A''$ and it is an exercise to show that for codimension 1, this situation does happen in an essentially unique way, and that it does not happen for $k=2$ (reason: such subspaces are affinely coindependent and so an associated linear system expressing their intersection has a solution).

Assume that neither $A$ nor $B$ is a midset.  In case both are long, we may replace with complements to assume both are short.  In any case, we may assume that $A$ is short, of cardinality $2^{d-1}-2^{d-k-1}$, for some integer $k$, $0<k\le \frac d 2$.  

First assume that $B$ is short, say of cardinality $2^{d-1}-2^{d-r-1}$, for $r>k$.  
Then $A+B$ has cardinality 
$2^{d-k-1}-2^{d-r-1}=2^{d-r-1}(2^{r-k}-1)$.  
Since $A+B$ is short, 
there exists an integer $s\le \frac d 2$ so that 
$2^{d-r-1}(2^{r-k}-1)=2^{d-1}-2^{d-s-1}=2^{d-s-1}(2^s-1)$.  
If both sides are powers of 2, then 
$r=k+1$, $s=1$ and $d-r-1=d-s-1$ implies that $r=s=1$ and $k=0$, a contradiction.
Therefore both sides are not powers of 2 and so $r=s$ and $s=r-k$  and so $s=r$ and $k=0$, a final contradiction.

Therefore  $B$ is long, of cardinality 
$2^{d-1}+2^{d-r-1}$, for $r>0$.  
Then $A+B$ has cardinality $2^{d-r-1}+2^{d-k-1}$.  Since $r\ge 1, k\ge 1$, this number is at most $2^{d-1}$ and is less than $2^{d-1}$ if $(r,k)\ne (1,1)$.  

Suppose that $r=k$.  Then $2^{d-r-1}+2^{d-k-1}=2^{d-r}$ is $2^{d-1}$ or $2^{d-2}$, implying $r=k=1, r=k=2$, respectively.  

Suppose that $r<k$.   Then  $A+B$ is short and 
there exists an integer $s\le \frac d 2$ so that $2^{d-r-1}+2^{d-k-1}=2^{d-1}-2^{d-s-1}$, and 
$2^{d-k-1}(2^{k-r}+1)=2^{d-s-1}(2^s-1)$.  Now, 
$2^{k-r}+1$ is odd, so it follows that $s=k$, $k-r=1$ and $s=2$.  So, $k=2, r=1$.

Suppose that $r>k$. Then $A+B$ is short and 
there exists an integer $s\le \frac d 2$ so that $2^{d-r-1}+2^{d-k-1}=2^{d-1}-2^{d-s-1}$, and 
$2^{d-r-1}(2^{r-k}+1)=2^{d-s-1}(2^s-1)$.  It follows that $s=r$, $r-k=1$ and $s=2$.  So, $k=1, r=2$.
\eop

\subsection{About defect 1 midsets} 

\begin{lem}\labtt{midsetmeethyp}  Let $d\ge 3$.  
Suppose that $B$ is a midset of defect 1.  Then $B$ contains affine hyperplanes of codimension 2.  Suppose that $A$ is an affine codimension 2 space contained in $B$.  There exists a unique hyperplane $H$ so that $B\cap H=A$.  (The other two hyperplanes which contain $A$ are cleansing hyperplanes for $B$ \reftt{cleansing}.) 
\end{lem} 
\pf
Let $A$ and $A'$ be any pair of disjoint codimension 2 subspaces.  
Then $A+A'$ is a midset and it has defect 0, 1 or 2 if $A'$ has a translate  which meets $A$ in codimension 0, 1 or 2, respectively.  The first statement follows from \reftt{cosetcleanliness} and transitivity of $AGL(d,2)$ on midsets of a given defect \reftt{localdirty}.  

For the second, consider the three hyperplanes $H_1, H_2, H_2$ which contain $A$.   
Suppose that $H_1\cap B > A$.  
Then $|H_1+B|=|H_1|+|B|-2|H_1\cap B|=2^d-2|H_1\cap B|<2^{d-1}$, whence $H_1$ is a cleansing hyperplane, and so $|H_1+B|=2^{d-1}-2^{d-1-1}=2^{d-2}$ and $|H_1\cap B|=\half (2^{d-1}+2^{d-1}-\cdot 2^{d-2})=2^{d-3}3$.  This means that at most two of the $H_i$ meet $B$ in a set larger than $A$.  Therefore, 
since $H_i\setminus A$ for $i=1, 2, 3$, partition $\FF_2^d \setminus A$, exactly two of the 
$H_i$ meet $B$ in a set larger than $A$ and so 
there exists an $H$ which meets $B$ in $A$, and by above counting, it is unique.  \eop

\section{More group theory for BRW groups}

We 
list some assumed results from group theory.

\begin{lem}\labtt{mindimsp} 
(i) A faithful module for 
$\prod _{1}^k Sym_3$ 
in characteristic 2 has dimension at least $2k$.  

(ii) 
A faithful module for $Sp(2k,2)$ in characteristic 2 has dimension at least $2k$.  
\end{lem}
\pf  Let $K_1\times \dots \times K_k$ be the natural direct product of $K_i\cong Sp(2,2)\cong Sym_3$ in $Sp(2k,2)$.  Clearly, (i) implies (ii).  We prove (i).

We may assume that the field $F$ is algebraically closed and that $k\ge 2$.  
Let $M$ be a module of minimal dimension.  
 Consider the decomposition $M=M'\oplus M''$, where $M'=[M,O_3(K_1)]$ and $M''=C_M(O_3(K_1))$.  

Clearly, $dim(M')$ is a positive even integer.  Suppose $M''\ne 0$.  Then by induction applied to the action of 
$K_2 \times \dots \times K_k$ on $M''$, we have $dim(M'')\ge 2(k-1)$ and we are finished.  Suppose $M''=0$.  
Then we may decompose $M''=P\oplus Q$ where $P$ and $Q$ represent the two distinct linear characters of $O_3(K_1)$.  The actions of $K_2 \times \dots \times K_k$ on $P$ and $Q$ are faithful and equivalent since $P$ and $Q$ are interchanged by elements of $K_1$.  We now finish by induction.  
\eop

\begin{lem}\labtt{savep} 
Let $\FF_2^{2m}$ have a nonsingular quadratic form 
of type $\nu = \pm$ and let $sv(m,\nu)$, $av(m, \nu)$ denote the number of singular and nonsingular vectors in the case of type $\nu = \pm$.  Then 
$sv(m,\nu)=(2^m-\nu 1)(2^{m-1}+\nu 1)$ and  
$av(m, \nu)=(2^m-\nu 1)2^{m-1}$. 
\end{lem}
\pf Well-known.  Note that $sv(m,\nu)+av(m, \nu)+1=2^{2m}$.  \eop

\begin{lem}\labtt{orbits} Let $k\ge 2$.  
Let $U$ be the essentially unique $2k+1$ dimensional 
$\FF_2$-module for $Sp(2k,2)$ with socle of dimension 1 and quotient the natural $2k$-dimensional module.  Then 
(i) $U$ is the natural module for $O(2k+1,2)$; 
(ii) The orbits of $Sp(2k,2)$ on $U$ consist of the two 1-point orbits lying in the radical, and the singular points and the nonsingular points.  Each of the latter orbits form coset representatives for the nontrival cosets of the radical.
\end{lem}
\pf  This is mainly the 1-cohomology result \cite{Poll}, plus a standard interpretation of $Ext^1$.  
\eop

\subsection{For clean involutions}  

We use the following notation throughout this subsection. 

\begin{nota}\labttr{qt}
We have the clean upper  involution $t$ of defect $k \ge 1$. 
 Take a  CMZ decompostion $C_R(t)=PZ$.   
Denote by $q_t$ the quadratic form on $Z=Z_t$ described in \reftt{orbits}.  The subscript indicates dependence on the involution, $t$.  Call $z\in Z$ {\it singular} or {\it nonsingular}, according to the value of $q_t(z)$.  
\end{nota}

\begin{lem}\labtt{zindecmodule}  Use the notation of \reftt{qt}.  
For all $k\ge 1$, the set map $x\mapsto [x,t]$ takes $R\setminus C_R(t)$ to the set of nonsingular vectors in $Z$ with respect to the invariant quadratic form.  

For $k\ge 2$, the action of $C_G(t)$ as 
$Sp(2k,2)$ on $Z$ is indecomposable; the  upper L\"owey series has factors of dimensions 1, $2k$.  
\end{lem}
\pf
Let $f$ be the commutator map $R\rightarrow Z$ defined by $f(x):=[x,t]$.  Every coset of $Z(R)$ in $Z$ contains an element of $Im(f)$. 
If $f(x)=f(y)$, we have $1=f(x)f(y)=[x,t][y,t]$, which is congruent to $f(xy)$ modulo $\la -1 \ra$. 
If $f(xy) \in \la -1 \ra$, then $xy \in C_R(t)$.  Therefore $f$ maps $R/C_R(t)$ isomorphically onto $Z/\la -1 \ra$.    
Also, the image of $f$ is a set of cardinality $2^{2k}$ which contains $1$ and is invariant under $C_R(t)$, which  acts on $Z$ as $Sp(2k,2)$, i.e., $Im(f)$ is a $C_R(t)$-invariant transversal to $Z(R)$ in $R$.  

We compute that (*) $f(xy)=[xy,t]=[x,t]^y[y,t]=f(x)^yf(y)$.  

We claim that $Z$ is an indecomposable module for $Sp(2k,2)$.  Suppose it is decomposable.  
Then $Im(f)$ must be either a subspace of $Z$ complementing $Z(R)$ or essentially a coset of some $C_R(t)$-invariant subspace, say $Z_0$, namely it is the set $Y$ which is the notrivial coset with $-1$ replaced by 1.    
Then there exists a homomorphism $h: R\rightarrow Z_0$ with the property that $f(x)=-h(x)$ if $x \not \in C_R(t)$ and $h(x)=1$ if $x\in C_R(t)$.   

It can not be a subspace since $[R,t]$ is normal in $R$.  So, the second alternative applies to $Im(f)$.  Now, we shall get a contradiction, using (*).  

Note that we have an alternating  bilinear form $g$ on $Z$ with values in $Z(R)$, defined by $g(a,b):=[a', t, b']$ where priming on $a\in Z$ means an element $a'\in R$ so that $f(a')=a$.  It helps to think of the Hall commutator identity $[x,y^{-1},z]^y[y,z^{-1},x]^z[z,x^{-1},y]^x=1$.  

There is a $g$-totally singular subspace of dimension $k+1$ in $Z$, say $W$. 
Assuming that $Im(f)=Y$, we take any elements $a, b, c$ in $R$ so that $abc=1$ and none of $a, b, c$ is in $C_R(t)$.  
Then $f(a)f(b)f(c)=(-1)^3h(a)h(b)h(c)=-1$.  From (*), we get $f(c)=f(ab)=f(a)^bf(c)$.  Now choose $a, b, c$ so that 
$f(a), f(b), f(c) \in W$ (this is possible since $k \ge 2$).  Then $g(f(a),f(b))=1$ implies that 
$f(a)^b=f(a)$, which implies that $f(c)=f(a)f(b)$, in contradiction with $f(a)f(b)f(c)=-1$. This proves that  $Z$ is indecomposable.

At this point, we know that $Im(f)$ is one of two orbits for $Sp(2k,2)$ in $Z$, the singular one and the nonsingular one.  We claim that it is the singular one.  Suppose otherwise.  Take $W$ and $a, b, c$ as above.  Then (*) implies that (in additive notation) the sum of two orthogonal nonsingular vectors is nonsingular, a contradiction. 
\eop 

\begin{de}\labttr{typephi}  Let $\phi$ be a linear character of $Z$ which is nontrivial on $Z(R)$.  Then $Ker(\phi )$ is 
a nonsingular quadratic space by restriction of $q_t$ \reftt{qt}.  Its {\it type} is plus or minus, according to the Witt index of the restriction of $q_t$.  
\end{de}

\begin{lem}\labtt{y} 
Consider 
$X:=\{(\varphi , z) | \varphi \in Hom(Z,\FF_2 ), z \in Z\}$ and let 
$Y_{\vep, \zeta, \eta} :=\{(\varphi , z)\in X  | \varphi (Z(R))\ne 1, z\ne 1, q_t(z)=\zeta, type(\varphi )=\vep,  \varphi (z) = \eta \}$, for $\zeta, \eta \in \FF_2$.  
Then $C(t)$ is transitive on $Y_{\vep, \zeta, \eta}$, for $\zeta, \eta = 0, 1$.  

The orbit lengths are 

$|Y_{\vep, 0, 0} |=(2^{2k-1}+\vep 2^{k-1})sv(k,\vep )$;   

$|Y_{\vep, 0, 1} |=(2^{2k-1}+\vep 2^{k-1})av(k,\vep )$;

$|Y_{\vep, 1, 0} |=(2^{2k-1}+\vep 2^{k-1})av(k,\vep )$;   

$|Y_{\vep, 1, 1} |=(2^{2k-1}+\vep 2^{k-1})sv(k,\vep )$.  

Note that rows 2 and 3 are equal and rows 1 and 4 are equal.  

\end{lem}
\pf
It is well-known that $C(t)/O_2(C(t))\cong Sp(2k,2)$ acts with two orbits on characters of $Z$ 
which take nontrivial value on $Z(R)$.  These orbits have respective stabilizers the natural subgroups  $O^\vep (2k,2)$ and respective lengths $2^{2k-1}+\vep 2^{k-1}$.  
The rest follows from \reftt{savep}.  
\eop

\begin{nota}\labttr{gvep} Let $t\in G$ be an involution.  Then $C_G(t)$ acts on each eigenlattice $\lvept$.  Its image in $O(\lvept)$ is denoted $G_\vep$. 
\end{nota}

\begin{lem}\labttr{irrcgt} The action of $C_G(t)$ on $L^\pm (t)$ is irreducible.  The center of $G_\vep$ is just $\{\pm 1\}$.
\end{lem}
\pf  The second statement follows from orthogonality of the representation plus absolute irreducibility, which we now prove.  
We prove irreduciblity for a natural subgroup of $C_G(t)$ of the form 
$AB$, where $[A,B]=1$, $A \cong 2^{1+2(d-k)}$, $Z\le B,  B/Z\cong Sp(2k,2)$; see \reftt{cmz}, \reftt{zindecmodule}.  
Every faithful irreducible of $A$ has dimension $2^{d-2k}$.  The central involution of $R$ is in $Z$ and so every irreducible of $B$ on $\QQ \otimes L$ involves an orbit of characters of $Z$ of cardinality $2^{2k-1}\pm 2^{k-1}$, and both orbit lengths occur with multiplicity $2^{d-2k}$.  Therefore, just two irreducibles for $AB$ occur in $\QQ \otimes L$, and they have respective dimensions  
$2^{d-2k}(2^{2k-1}\pm 2^{k-1})=2^{d-1}\pm 2^{d-k-1}$.  The conclusion follows.  
\eop

\begin{lem}\labtt{tracez} Assume $t$ is clean with positive trace.  
Let $z \in Z$, $z\ne \pm 1$.  The trace of $z$ on $L^\pm (t)$ is 
$\pm 2^{d-k-1}$ if $z$ is $q_t$-singular and is 
$\mp 2^{d-k-1}$ 
if $z$ is $q_t$-nonsingular.  
\end{lem}
\pf 
We use the subgroup denoted $AB$ in the proof of \reftt{irrcgt}.  For $AB$, the module $\lvept$ decomposes as a tensor product of irreducibles.  It suffices to prove that the trace of $z$ on the tensor factor irreducible for $B$ is $\pm 2^{k-1}, \mp 2^{k-1}$, respectively.  

Note that $2^{2k}-1-sv(k,\vep )=(2^k-\vep) (2^k - \vep -(2^{k-1}+\vep ))=(2^{k}-\vep )2^{k-1}$.

We use \reftt{savep} to deduce that 
$$|Y_{\vep, 0, 0}|=2^{k-1}(2^k+\vep)sv(k,\vep)=2^{k-1}(2^k+\vep)(2^k-\vep )(2^{k-1}+\vep )$$
and 
$$|Y_{\vep, 0, 1}|=2^{k-1}(2^k+\vep )(2^{2k}-1-sv(k,\vep))=2^{k-1}(2^k+\vep )(2^k-\vep)2^{k-1}.$$

Let $\Phi_\vep :=\{ \varphi | \varphi (Z(R)) \ne \{1\} \}$.  
A  given singular $z\in Z$ is in the kernel of 
$2^{k-1}(2^{k-1}+\vep )=2^{2k-2}+\vep 2^{k-1}$ characters in  $\Phi_\vep$ and outside the kernel of 
$2^{2k-2}$ characters in  $\Phi_\vep$.  It follows that 
the trace of $z$ on $L^\vep (t)$ is $\vep 2^{k-1}$.  
Singular and nonsingular elements of $Z\setminus Z(R)$ are paired by congruence modulo $Z(R)$.  Therefore, nonsingular elements have trace $-\vep 2^{k-1}$.  
\eop  

\subsection{For dirty involutions}

We assume the following notation throughout this subsection.  

\begin{nota}\labttr{ulfactorization} Let $t$ be a dirty split upper involution of defect $k$.  A {\it UL factorization of $t$} is an expression $t=u\ell$, where $u$ is a clean involution and $\ell$ is a lower involution (note that all of $t, u, \ell$ commute).  Write ${\cal UL}(t)$ for all pairs $(u,\ell )$ as above.  Let ${\cal U}(t)$ be the set of $u$ and let ${\cal L}(t)$ be the set of $\ell$ which arise this way.  We have $|\{{\cal UL}(t)\}|=2^{{1+2(d-2k)}+2k}-2^{1+2k}$.  
\end{nota} 

We get a result for traces of $u$ and $\ell$ on $\lvept$ which is similar to \reftt{tracez}.  

\begin{lem}\labtt{ydirty} 
On $\lvept$, the trace of $z\in Z\setminus Z(R)$ is 0 and the trace of $\ell$ is $\pm 2^{d-k-1}$, for all $\ell \in {\cal L}(t)$.  
\end{lem} 
\pf
We assume $\vep=+$ (the other case is similar). 
It suffices to consider the sublattices $L(a,b)$, where $u$ acts as $a$ and $\ell$ acts as $b$.  Recall that the eigenlattices for $\ell$ are $ss\bw {d-1}$ lattices, for which we may use \reftt{y} to compute the traces for $z$.   Without loss, we may assume that $z$ has nonnegative traces.       We get:  
$$\begin{matrix}
sublattice & rank & multiplicity \ of &  multiplicity \ of  \cr
& &  +1 \ for \ z &  \ -1 \ for \ z \cr
L(+1,+1)  & 2^{d-2}+2^{d-k-2} &  2^{d-3}+2^{d-k-2} & 2^{d-3} \cr
L(+1,-1)   & 2^{d-2}+2^{d-k-2} &  2^{d-3}+2^{d-k-2} & 2^{d-3} \cr
L(-1,+1)   &  2^{d-2}-2^{d-k-2} & 2^{d-3}-2^{d-k-2} & 2^{d-3} \cr
L(-1,-1)   &  2^{d-2}-2^{d-k-2} & 2^{d-3}-2^{d-k-2} & 2^{d-3} \cr
\end{matrix}$$
\eop 

\section{About inherted groups} 

We continue to use the notations 
$G:=\gd d$, $R:=\rd d$.  See  the  ancestor section of 
\cite{bwy} for discussion.

\begin{nota}\labttr{bargeps}  
We use  bars for images under restriction $C_G(t)\rightarrow O(\lvept )$.  
As in \reftt{gvep}, we write $G_\vep$ for the image of $C_G(t)$ in $O(\lvept )$  under the restriction homomorphism. 
\end{nota}

\begin{lem}\labtt{ancestralradical}
Suppose that $Z$ is an elementary abelian subgroup of $R$ containing $Z(R)$ and that $rank(Z)=s+1$.  Let $L_\l$ be the eigenlattice for $L$, defined by the linear character $\l$ of $Z$, which is assumed to be nontrivial on $Z(R)$.  The set $\cal F$ of such $\l$ has cardinality $2^s$.  

There is a  finite subgroup of the orthogonal  group $O(\QQ \otimes \lvept )$ 
of the form 
$\prod_{\l} R_\l$ with the property that $R_\l$ acts on $L_\l$ as a lower group and acts trivially on $L_\mu$ for $\mu \ne \l$.  We have $|\prod_{\l} R_\l |=2^{s(1+2(d-s))}$.

(i) When $s = 1$, $\overline{C_G(Z)}\ge \prod_{\l} R_\l$.

(ii) When $s = 2$, $\overline{C_G(Z)}\cap \prod_{\l} R_\l$ is an index $2^{2(d-2)}$ subgroup of $\prod_{\l} R_\l$ with the property that if $\cal J$ is any 3-set in $\cal F$, then the projection of 
$\overline{C_G(Z)}$ to 
$\prod_{\l \in {\cal J}} R_\l$ is onto.  The kernel of this homomorphism is just $Z(R_\mu )$, where $\mu \in {\cal F}$ is the index missing from $\cal J$.  
\end{lem}

\begin{lem}\labtt{inradical} 
Let ${\cal I} \subseteq {\cal F}$ be any nonempty collection of characters as in \reftt{ancestralradical} and let $J:=J({\cal I})$ be the direct summand of $L$ determined by $span\{J_\nu | \nu \in {\cal I}\}$.  If $\l \in {\cal I}$ and $g \in C_G(Z)$ acts trivially on $J_\l$, then $g$ acts on $J$ as an element of the group $\prod_{\l} R_\l$, defined in \reftt{ancestralradical}.  
\end{lem}
\pf 
If $\mu, \nu$ are any two distinct indices so that $L_\mu$ and $L_\nu$ are stable under $h \in G$, then if $h$ acts trivially on $L_\mu$ modulo its first lower twist, then $h$ does the same on $L_\nu$.  By considering  all distinct pairs of indices $\mu, \nu \in {\cal I}$, we deduce that $g$ acts on $J$ as a member of $\prod_{\eta \in {\cal I}} R_\l$.  See \cite{bwy} 
\eop  

\begin{coro}\labtt{normzinherited} 
Use the notation of \reftt{inradical}.  Assume that $s = 2$, $\cal I$ has cardinality 3 and 
$N_{O(J)}(\overline Z)=\overline{N_G(Z)}C_{O(J)}(\overline Z)$.  Then 
$N_{O(J)}(\overline Z)$ is inherited.  
\end{coro}
\pf \reftt{inradical} and \reftt{ancestralradical}(ii).  \eop

\begin{coro}\labtt{smalldefectinherited} 
Let $t\in G$ be an involution and let $Z:=Z(C_R(t))$.  
(i) 
Suppose that $t$ is a clean involution of defect $1$.  Then 
$N_{O(\lvept)}(\overline Z)$ is inherited.

(ii) 
Suppose that $t$ is a split dirty involution of defect 1.  Then 
$N_{O(\lvept)}(\overline Z)$ is inherited.
\end{coro} 
\pf 
Note that defect 1 implies that $s=2$, in the notation of \reftt{ancestralradical}.  
(i): 
This  follows from \reftt{normzinherited}.

(ii): 
Let $t$ be such an involution. Let $t=u\ell$ be a UL-factorization \reftt{ulfactorization}.  
We define $Z_u$ as $Z(C_R(u))$ and define $Z:=Z(C_R(t))=Z_u\times \la \ell \ra$.  A character value analysis shows that elements of $Z_u$ have 0 trace on $\lvept$ and elements of the coset $Z_u\ell$ have nonzero trace \reftt{ydirty}. 
Therefore, $N_{O(\lvept)}(\overline Z_u) \ge N_{O(\lvept)}(\overline Z)$.  

We shall use \reftt{normzinherited} to prove that $N_{O(\lvept)}(\overline Z_u)$ is inherited by showing that the latter group induces only $Sp(2,2)$ on $Z$.  
Assume that this is 
false.  We have an action of $AGL(2,2)\cong Sym_4$ on $Z$.  Let $H$ be the linear group which 
$N_{O(\lvept)}(\overline Z_u)$ induces on $Z_u$.  The action of $N_{C_G(t)}(Z_u)$ on $Z_u$ preserves the coset $Z_u\setminus Z$ and has orbits modulo $Z(R)$ of lengths 1 and 3.  Its orbits on 
 $Z_u\setminus Z$ must have lengths 1,1,3,3 since elements in that coset have nonzero trace on $\lvept$ so are not conjugate to their negatives.  It follows that a Sylow 2-group $S$ of 
 $N_{O(\lvept)}(\overline Z_u)$ fixes an element, say $\ell$,  in this coset.  If $x \in Z \setminus Z(R)$, then there exists $g \in S$ so that $x^g=-x$, since we are assuming an action of $AGL(2,2)$ on $Z$.  It follows that $(x\ell )^g=-x\ell$, which is a contradiction since $(x\ell, xu)$ is a UL factorization of $t$ (because $xu\in uZ=u^R$ consists of clean elements).  
\eop

\section{The split defect 1 cases}

\subsection{The clean defect 1 case}

\begin{de}\labttr{ssdsubgroup}  Suppose that $M$ is an integral lattice and $X$ is a SSD lattice.  Define $SSD(M,X)$ to be the subgroup of $O(M)$ generated by the SSD involutions associated to sublattices of $M$ which are isometric to $X$. 

This is clearly a normal subgroup of $O(M)$.
\end{de}

 We continue to use the notation \reftt{cmz}.  Since the defect is 1, $rank(Z)=3$.  

\begin{rem}\labttr{xperp}  In the notation of \reftt{ssdsubgroup}, if $X$ is SSD and $det(M)=1$, then $M \cap X^\perp$ is SSD.  This will apply for us when $M\cong \bw d$ and $d$ is odd.  
\end{rem}  

\begin{lem}\labttr{cleandefect1z}  Suppose that $d>3$.  
Let $t$ be a clean involution of defect 1 and positive trace.  Then $SSD(\lpt, ss\bw {d-1})=Z$.   
\end{lem}  
\pf 
A sublattice $X$ of $\lpt$ which is isometric to $ss\bw {d-1}$ is SSD in the overlattice $L$.  By \cite{bwy}, the associated SSD involution is lower (here, we are using $d>3$), so lies in $C_R(t)$.  Since $Tr_{\lpt}(\vep_X) \ne 0$ (see \reftt{tracez}), $\vep_X \in Z$, the only elements of $C_R(t)$ which have nonzero trace on $\lpt$.  The action of $C_G(t)$ on $Z$ is that of $O(2k+1,2)$ on its natural module \reftt{zindecmodule}.  Therefore, every element of $Z\setminus Z(R)$ is such an SSD involution.  
\eop 

\begin{lem}\labttr{cleandefect1inh} 
  Suppose that $d>3$.  
Let $t$ be a clean involution of defect 1 and positive trace.  Then $O(\lpt)$ is inherited.  
\end{lem} 
\pf
By \reftt{cleandefect1z}, $\overline Z$ is normal in $O(\lpt)$.  Now use \reftt{normzinherited}
\eop 

\begin{rem}\labttr{cleandefect1postrace}  
If $t$ is a clean involution of defect 1 and positive trace, $\lmt \cong ss\bw {d-1}$, whose automorphism group is known. 
\end{rem}

\subsection{The split dirty defect 1 case}

\begin{lem}\labtt{dirtydefect1z}  Suppose that $d>3$ and $d$ is odd.  
Let $t$ be a split dirty involution of defect 1.   Then 
$SSD(\lvept, ss\bw {d-2})= \overline{C_R(t)}$ and $\overline Z$ is a subgroup of $Z(SSD(\lvept, ss\bw {d-2}))$ which is normal in $O(\lvept)$.  
\end{lem} 
\pf 
We may suppose that $t=\vep_b$ for a defect 1 midset $b\in RM(2,d)$ and we may assume that $\vep =+$.   
Since $d$ is odd, a $ss\bw {d-2}$ sublattice is SSD.  
Let $X$ be such a sublattice of $\lvept$.  Its associated involution in $O(L)$ is conjugate to an involution of the form $\vep_c \in {\cal E}$, where the codeword $c$ is an affine codimension 2 subspace.

Since $\vep_c$ acts nontrivially on $\lpt$,  $c\cap  b=\emptyset$   Let $b'$ be the complement of $b$.  
We may consider the involution $\vep_h$ where $h$ is a hyperplane so that $h \cap b'=c$ \reftt{midsetmeethyp}.  Then $\vep_c$ acts on $\lpt$ as $\vep_h$, which is a lower involution.

Define $K$ to be the normal subgroup of $O(\lpt)$ generated by all $\overline{\vep_X}$, where $X$ is a SSD sublattice isometric to $ss\bw {d-2}$.  This is a subgroup of $\overline{C_R(t)}$ which is normal in 
$\overline{C_G(t)}$ and contains $\overline{\vep_h}$, so is not contained in $\overline{Z(R)}$.  The normal subgroups are $\overline 1, \overline{Z(R)}, \overline Z,  \la \overline Z, \overline {\ell}  \ra, \overline{C_R(t)}$, where $\ell \in {\cal L}(t)$ is any lower part of a UL-factorization.  For all such normal subgroups, $Y$ not contained in $Z(R)$, we claim that $Z$ is normal in $N_{O(\lpt )}(Y)$.  This is obvious except when $Y$ is one of the latter two cases.  In those cases, $Z(Y)=\la Z, \ell \ra$.  In the action on $\lpt$, the elements of $Z\setminus Z(R)$ have trace 0 and the elements of $Z\ell$ have nonzero trace (see \reftt{ydirty}).  The claim follows and so does the lemma since $K$ is normal in $O(\lpt )$.  
\eop

\begin{lem}\labtt{dirtydefect1inh}  Suppose that $d>3$ and $d$ is odd.  
Let $t$ be a split dirty involution of defect 1.   Then $O(\lvept)$ is inherited.  
\end{lem} 
\pf 
Use 
\reftt{dirtydefect1z} 
and 
\reftt{smalldefectinherited}(ii).  
\eop

This completes the proof of \reftt{maintheorem5}.  

\section{The nonsplit defect 1 case}

The style of proof here is rather different.  The smallest value of $d$ for this case is $d=2$.  Involutions in $\brw 2 \cong W_{F_4}$ are discussed in \reftt{alldefectsoccur}.  

\begin{lem}\labtt{nonsplitd=3} 
Let $t$ be a nonsplit involution of defect 1 in $\brw 3$.  Then $\lvept \cong L_{A_1^4}$.  
\end{lem}
\pf Let $L=\bw 3 \cong \leh$.

{\it First proof: }
In a root system of type $E_8$, there are $A_1^8$ subsystems.  Let $\Psi$ be one.  If $F$ is a linearly independent 4-set in $\Psi$, either $L[F]:=\QQ\otimes  F\cap L\cong L_{D_4}$ or $L_{A_1^4}$, and both occur.  Let $F$ be the latter type.  The sublattice $L[F]$ is SSD, and in our accounting so far, has not appeared as an eigenlattice of an involution on $\bw 3$.  It therefore must represent the missing type.  

{\it Second proof: }  There exists a lower   dihedral group $D\le C_R(t)$.  Let $u, v$ be involutions which generate $D$.  Then by 2/4-generation \cite{bwy}, $L=L^\pm(u)\oplus L^\pm(v)$, all summands are  $\bw {2} \cong L_{D_4}$ lattices which are $t$-invariant and on them $t$ acts like a nonsplit dirty involution.   Each eigenlattice has an orthogonal basis of norms 2, 4 (see \reftt{alldefectsoccur}) and the total eigenlattice of the restriction of $t$ has index 4.  It follows that  $|L:Tel(L,t)|=|L^\vep (u):Tel(L^\vep (u),t)|^2=4^2$ whence $det(Tel(L,t))=2^8det(L)=2^8$ and $det(\lvept )=2^4$.    From the section on actions of 2-groups in \cite{bwy}, we conclude that $Tel(L,t)=[L,t]+2L$, whence each $\lvept$ is isometric to $\sqrt 2 P$, where $P$ is an integral lattice of determinant 1 and rank 4.  Therefore $P$ is the square unimodular lattice by the well-known classification of integral unimodular lattices of rank at most 8.  This completes the proof.  
\eop 

\begin{lem}\labtt{jcfnonsplitk=1} If $t$ is a nonsplit involution of defect 1, $L/2L$ is a free $\FF_2\la t \ra$-module, i.e., the Jordan canonical form for $t$ consists of $2^{d-1}$ blocks of degree 2.
\end{lem}
\pf The result may be checked directly for $d\le 2$ since 
we know $Tel(t)$ and $Tel(t)+2L/2L$ is the fixed point space for the action of $t$ on $L/2L$ (see \cite{bwy}).  
The idea is to use induction on $d$ plus the fact that $t$ leaves invariant the summands of a decomposition $L=L^\pm(u)\oplus L^\pm(v)$, where $u, v$ generate a lower dihedral group which centralizes $t$.  This proves that $L$ is a free $\ZZ \la t \ra$-module, so reduction modulo 2 has the claimed structure.  
\eop

\begin{lem}\labtt{nonsplitevenness} Let $t$ be a nonsplit involution of defect 1.  Then $\lvept$ is doubly even for $d \ge 4$, i.e., $\frac 1 {\sqrt 2}\lvept$ is an even integral lattice.   
\end{lem}
\pf 
When $d=4$,  $L^\pm (u) \cong \rtleh$, so the property is clearly true. 
By \reftt{jcfnonsplitk=1}, we have $Tel(t)=2L+[L,t]$.  For $x, y \in L$, $(x(t-1), y(t-1))=(x,y)+(xt,yt)-(x,yt)-(xt,y)=2(x,y)-2(x,yt) \in 2\ZZ$.  
It follows that $(Tel(t),Tel(t))\le 2\ZZ$.  We take $x=y$.  We want $(x,xt)\in 2\ZZ$ to conclude that $x(t-1)$ has norm divisible by 4.  This will follow if it is so for a spanning set.  
Consider the summands of a decomposition $L=L^\pm(u)\oplus L^\pm(v)$, where $u, v$ generate a lower dihedral group which centralizes $t$.  For $x\in L^\pm (u)$, which is a $ss\bw {d-1}$,  $x(t-1)$ has norm divisible by 4 for $d \ge 5$, by induction.  
\eop

\begin{lem}\labtt{nonsplitd=4} 
Let $t$ be a nonsplit involution of defect 1 in $\brw 4$.  Then $\lvept \cong \sqrt 2 \bw 2 \perp  \sqrt 2 \bw 2 \cong  \sqrt 2 L_{D_4} \perp \sqrt 2 L_{D_4}$.  
\end{lem}
\pf Let $L=\bw 4 \cong \leh$.
We follow the strategy in the 
proof of \reftt{nonsplitd=3}.  
 Then there exists a lower dihedral group $D\le C_R(t)$.  Let $u, v$ be involutions which generate $D$.  
 Then by 2/4-generation \cite{bwy}, 
 $L=L^\pm(u)\oplus L^\pm(v)$, all summands are  
$ss\bw {3}\cong \rtleh$ lattices which are $t$-invariant and on them $t$ acts like a nonsplit dirty involution. 
It follows that each $L^\pm(w)^\vep (t)$ is isometric to $\sqrt 2 L_{A_1^4}$, for any noncentral involution $w\in D$. Reasoning as in \reftt{nonsplitd=3}, we argue that $Tel(t)$ has index $2^8$ in $L$ and $det(Tel(t))=2^{16}det(L)=2^{24}$.  
From \reftt{nonsplitevenness}, we know that $Tel(t)$ is doubly even.  Therefore, each $\lvept$ is doubly even and has determinant $2^{12}$.  Therefore, there is an even integral lattice, $P$, so that $\lvept \cong \sqrt 2 P$,  $det(P)=2^4$ and $P$ contains a sublattice $Q$ isometric to $L_{A_1^8}$, of index $4$ in $P$.  

Let $r_1, \dots , r_8$ be an orthogonal basis of roots for $Q$.  
Any nontrivial coset of $Q$ in $P$ consists of even norm vectors, so contains an element of shape $\half \sum_{i \in I} r_i$, where $I\subseteq \{1, 2, 3, 4, 5, 6, 7, 8\}$ and $|I|=4$ or 8 (note that $exp(P/Q)\ne 4$ since vectors of shape $\sum_{i=1}^8 \pm \fourth r_i$ have norm 1).  

Let $I, I'$ be any two 4-sets which arise as above.  We claim that they are disjoint.  Assume otherwise.  Since $P$ is even, $I \cap I'$ is a 2-set.  Then $P$ is isometric to $L_{D_6}\perp L_{A_1} \perp L_{A_1}$, whence $C_R(t)$ fixes the unique indecomposable orthogonal summand isometric to $L_{D_6}$.  This is impossible since $C_R(t)$ 
contains a subgroup of shape $2^{1+4}_+$, whose faithful irreducibles have dimension divisible by 4.   We conclude that there exists a partition $J', J''$ of 
$\{1, 2, 3, 4, 5, 6, 7, 8\}$ so that $J'$ and $J''$ are 4-sets and $P=P'\perp P''$, where $P':=\{x \in P |supp(x) \subseteq J'\}$,  
$P'':=\{x \in P |supp(x) \subseteq J''\}$ and $P'\cong P''\cong \bw 2 \cong L_{D_4}$.  
\eop 

\begin{prop} \labtt{nonsplitstory}  For all $d\ge 2$, if $t\in \brw d$ is a nonsplit dirty involution, then $\lvept \cong ss\bw {d-2}\perp ss\bw {d-2}$.
\end{prop} 
\pf  If $d=2$, this is true by the discussion in \reftt{alldefectsoccur}.  
For $d=3, 4$, we use \reftt{nonsplitd=3}, \reftt{nonsplitd=4}.    

Let $d\ge5$.  Then there exists a lower   dihedral group $D\le C_R(t)$.  Let $u, v$ be involutions which generate $D$.  Then by 2/4-generation \cite{bwy}, $L=L^\pm(u)\oplus L^\pm(v)$, all summands are  ss$\bw {d-1}$ lattices which are $t$-invariant and on them $t$ acts like a nonsplit dirty involution.  By induction, we know the eigenlattices for $t$ on each.

Consider $L^+(u)\perp L^-(u)$.  The involution $v$ interchanges the summands and acts trivially on $L/L^+(u)\perp L^-(u)$.  The same is therefore true for 
the actions of $v$ on $L^+(u)^\vep (t) \perp L^-(u)^\vep (t)$ and 
$\lvept / L^+(u)^\vep (t) \perp L^-(u)^\vep (t)$.

Since $d-1\ge 4$, induction implies that each $L^\pm (u)^\vep (t)$  is the orthogonal sum of two orthogonally indecomposable lattices.  
Furthermore, if $S$ is one of these two indecomposable direct summands of $L^+(u)^\vep (t)$, we deduce  that the same is true for the actions  of $v$ on 
$S\perp S^v$ and on $\lvept \cap (\QQ \otimes (S\perp S^v))/S\perp S^v$.

We finish by quoting the uniqueness theorem \cite{bwy, bwycorr}, applied to the containment of $S\perp S^v$ in $\lvept \cap (\QQ \otimes (S\perp S^v))$, for each $S$.  Note that $t$ centralizes a natural $\brw {d-2}$-subgroup of $\brw d$ and that it stabilizes $S$ and $S^v$.  
\eop

The main result \reftt{maintheorem6} follows.

\section{Appendix: About BRW groups.   }

This is an updated and corrrected version of Appendix 2 from \cite{bwy}.

Basic theory of extraspecial groups extended upwards by their
outer automorphism group has been
developed in several places.  We shall use 
\cite{GrEx, GrMont, GrDemp, GrNW, Hup, BRW1, BRW2, B}.

\begin{nota}\labttr{brw}  Let $R\cong 2^{1+2d}_\varepsilon$ be an
extraspecial group  which is a subgroup of $GL(2^d,\FF)$, for a field
$\FF$ of characteristic 0. Let $N:=N_{GL(2^d,\FF)}(R)\cong \FF^\times.2^{2d}
O^\varepsilon (2d,2)$.  The {\it Bolt-Room-Wall group} is a subgroup of
this of the form $2^{1+2d}_\varepsilon.\Omega^\varepsilon (2d,2)$.  
If $d\ge 3$ or $d=2,
\varepsilon =-$,  $N'$ has this property.  For
the excluded parameters, we take a suitable subgroup of such a group for
larger $d$.  
We denote this group by $BRW^0(2^d,\vep)$ or $\dg d$.  It is uniquely
determined up to conjugacy in  $GL(2^d,\FF)$ 
by its isomorphism type if  $d\ge 3$ or $d=2,
\varepsilon =-$.  It is conjugate to a subgroup of $GL(2^d,\QQ)$ if
$\varepsilon = +$.  Let $R=\rd d$ denote $O_2( \gd d )$.  We call
$\rd d$ {\it the lower group} of $BRW^0(2^d,+)$ and call $\gd d/\rd d$ {\it
the upper group} of $BRW^0(2^d,+)$.

For $g\in N$, define 
$C_{R \ mod \ R'}(g):=\{ x \in R | [x,g]\in R'\}$, 
$B(g):=Z(C_{R \ mod \ R'}(g))$ and let $A(g)$ be some subgroup of 
$C_{R \ mod \ R'}(g)$ which contains $R'$ and complements $B(g)$ modulo
$R'$, i.e., 
$C_{R \ mod \ R'}(g)=A(g)B(g)$ and $A(g)\cap B(g)=R'$.  
Thus, $A(g)$ is extraspecial or cyclic of order 2.  
Define 
$c(d):=dim(C_{R/R'}(g))$, $a(g):=\half |A(g)/R'|$, $b(g):=\half
|B(g)/R'|$. 
 Then $c(d)=2a(d)+2b(d)$.  
\end{nota}

\begin{coro}\labttr{agld2}  Let $L$ be any $\ZZ$-lattice invariant under
$H:=BRW^0(2^d,+)$.   Then $H$ contains a subgroup $K\cong AGL(d,2)$
and $L$ has a linearly independent set of vectors $\{ x_i | i \in
\Omega\}$ so that there exists an identification of $\Omega$ with
$\FF_2^d$ which makes the $\ZZ$-span of 
$\{ x_i| i \in\Omega \}$ a  permutation module for $AGL(d,2)$ on
$\Omega$.  
\end{coro}
\pf
In $H$, let $E, F$ be  maximal elementary abelian subgroups and let $K$
be their common normalizer.  It satisfies $K/R \cong GL(d,2)$.  Now, let
$z$ generate $Z(R)$ and let $E_1$ complement $\la z \ra$ in $E$ and 
$F_1$
complement $\la z \ra$ in $F$.  The action of $K$ on the hyperplanes of
$E$ which  complement $Z(R)$ satisfies $N_K(E_1)F=K,  N_F(E_1)=Z(R)$. 
Now consider the action of $N_K(E_1)$ on the hyperplanes of $F$ which
complement
$Z(R)$.  We have that $K_1:=N_K(E_1)\cap N_K(F_1)$ covers
$N_K(E_1)/E$.  Therefore, $K_1/Z(R)\cong GL(d,2)$.  Let $K_0$ be the
subgroup of index 2 which acts trivially on the fixed points on $L$ of
$E_1$, a rank 1 lattice.  So, $K_0\cong GL(d,2)$.    Let $x$ be a basis
element of this fixed point lattice.  Then the semidirect product 
$F_1{:}K_0$ is isomorphic to $AGL(d,2)$ and $\{ x^g|g \in F_1\}$ is a
permutation basis of 
its $\ZZ$-span.  
\eop

\begin{de}\labttr{cleandirty}  We use the notation of \reftt{brw}.   An
element
$x \in N$ is  {\it dirty} if there exists $g$ so that 
$[x,g]=xz$, where $z$ is an element of  order 2 in the center.
If $g$ can be chosen to be of order 2, call $x$ {\it really dirty} or
{\it extra dirty}.  If $x$ is not dirty, call $x$ {\it clean}.
\end{de}

\begin{lem}\labtt{outer}  Let $\FF_2^{2d}$ be equipped with a
nondegenerate quadratic form with  maximal Witt index.  The set of maximal
totally singular subspaces has two orbits under $\Omega^+(2d,2)$ and these
are interchanged by the elements of $O^+(2d,2)$ outside 
 $\Omega^+(2d,2)$.  
\end{lem}
\pf
This is surely well known.  For a proof, see \cite{GrElAb}. 
\eop

\begin{de}\labttr{defect}  
An involution in  $\brw d$ has {\it defect $k$} if its commutator space on the Frattini factor of the lower group has dimension $2k$.  The defect is an integer in the range $[0,\frac d 2 ]$.   Note that an automorphism of $\rd d$ has even dimensional commutator space on $\rd d/Z(\rd d)$ if and only if it is even; see  \cite{GrMont}, \cite{GrElAb}.  
\end{de}

\begin{de}\labttr{splitnonsplit}  An involution in $\brw d$ is {\it split} if it centralizes a maximal elementary abelian subgroup of $\rd d$, and is otherwise {\it nonsplit}.
\end{de} 

\begin{nota}\labttr{eta} Write $R=D_1\dots D_d$ as a central product of dihedral groups, $D_i$ of order 8.  The involution $\a_{d,r}$ in 
$\Aut(\brw d)$, defined up to conjugacy, acts trivially on $d-r$ of the $D_i$ and performs an outer automorphism on the other $r$ of them. 
When $r=2k$ is even, $\a_{d,2k}$ is represented in $\brw d$ by an involution $\eta_{d,2k,+}$ (see \ref{nonsplitinvols}).  
In case $r=2k<d$, we define  an involution $\eta_{d,2k,-}:=\eta_{d,2k,+}z$, where $z$ is a noncentral involution in the above product of the $d-r$ elementwise fixed $D_i$.  
\end{nota}

\begin{thm}\labtt{traces}
We use the notation of \reftt{brw}, \reftt{cleandirty}.  Let $g \in N$. 
Then $Tr(g)=0$ if and only
if $g$ is dirty.  Assume now that $g$ is clean and has finite order. 
Then 
$Tr(g)=\pm 2^{a(g)+b(g)}\eta$, where $\eta$ is a root of unity.   
If $g\in BRW(d,+)$,
we may take $\eta = 1$. 
Furthermore, every coset of $R$ in $BRW(d,\varepsilon )$ contains a
clean element and if $g$ is clean, 
the set of clean elements in $Rg$ is 
just $g^R\cup -g^R$.  
\end{thm} 
\pf 
\cite{GrMont}.  
\eop

\begin{lem} \labtt{conjorthoginvols} Suppose that $t, u$ are involutions
in $\Omega^+(2d,2)$, for $d \ge 2$.  Suppose that their commutators on the
natural module $W:=\FF_2^{2d}$ are totally singular subspaces of the same
dimension, $e$.  Suppose that $e<d$ or that $e=d$ and that $[W,t]$ and
$[W,u]$ are in the same orbit under $\Omega^+(2d,2)$.   Then $t$ and $u$
are conjugate.  
\end{lem}
\pf 
Induction on $d$.
\eop

\begin{coro}\labtt{conjclean}  Suppose that $t, u$ are clean
involutions in
$H$ so that  $Tr(t)=Tr(u)\ne 0$.  Then $t$ and $u$ are conjugate in $\gd d$, if their common defect is less than $\frac d 2$.  If the defects are $\frac d 2$, then there are two classes.  
\end{coro}
\pf We may assume that $t,u$ are noncentral.  
These involutions are not lower and have the same dimension of fixed
points on $R/R'\cong \FF_2^{2d}$.  Let $T, U\le R$ be their respective
centralizers in $R$.  Since both $t, u$ are clean, $[R,t]$ and $[R,u]$
are elementary abelian subgroups of $T, U$, respectively.  
From \reftt{conjorthoginvols}, we deduce that $Rt$ and $Ru$ are conjugate
in $\gd d$ if their common defect is less than $\frac d 2$ and there are two possible conjugacy classes in case of common defect $\frac d 2$.  We may assume that $Rt=Ru$.  Now use \reftt{traces} to deduce
that $t$ is $R$-conjugate to $u$ or $-u$. The trace condition implies
that $t$ is conjugate to $u$.  
\eop 

\begin{rem}\labttr {splitting}  The extension 
$1\rightarrow \rd d \rightarrow \gd d \rightarrow \Omega^+(2d,2) \rightarrow
1$ 
is nonsplit for $d\ge 4$.  This was proved first in \cite{BRW2}, then
later in \cite{BE} and in \cite{GrEx} (for both kinds of
extraspecial groups, though with an error for
$d=3$; see
\cite{GrDemp} for a correction).  The article \cite{GrEx} gives a
sufficient condition for a subextension 
$1\rightarrow \rd d \rightarrow H \rightarrow H/\rd d \rightarrow
1$ to be split, and there are interesting applications, e.g. to the
centralizer of a 2-central involution in the Monster \cite{Gr72}.  A general
discussion of exceptional cohomology in simple group theory is in
\cite{GrNW}.  
\end{rem}

\begin{lem}\labtt{annihsing} Let $V=\FF_2^{2d}$ have a nonsingular quadratic form, $q$, of plus type.  Let $W$ be an isotropic subspace, $U:=W^\perp$.  Then every nontrivial coset of $U$ contains singular and nonsingular vectors if $d>1$. 
\end{lem}
\pf
Suppose that $v+U$ is a coset which consists entirely of either singular or nonsingular vectors.  Then for all $x,y \in v+U$, $q(x+y)=(x,y)+q(x)+q(y)=(x,y)$.  Take $a, b \in U$ so that $a+b=x+y$.  Then $(x,y)=q(a+b)=q(a)+q(b)+(a,b)$.  Also $(x+a,y+a)=(x,y)$ implies that $0=(x,a)+(a,y)=(x+y,a)=(a+b,a)=(a,b)$.  It follows that for any two elements $a, b$ of $U$, $(a,b)=0$.  Since $U$ is the annihilator of $W$, $U=W$.  Let $Z:=\{x \in W| q(x)=0\}$, a 
subspace of $W$ of codimension 0 or 1.  
Suppose $d>1$. Let $x \in V\setminus W$.  
If there is $z\in Z$ so that $(x,z)=1$, then $x$ and $x+z$ have different values under the quadratic form.  If this fails to be so, then $dim(Z)=0$, i.e., $d=2$ and $W$ contains nonsingular vectors.  Then $x$ annihilates a nonsingular vector, $w \in W$ and so 
 $x$ and $x+w$ have different values under the quadratic form. 
\eop

\begin{lem}\labtt{nonsplitinvols} Let $V=\FF_2^{2d}$ and let $g$ be an involution in $\Omega^+(2d,2)$ so that $[V,g]$ has dimension $r>1$ and contains nonsingular vectors.  There exists a basis of singular vectors $x_1, \dots , x_d, y_1, \dots  y_d$ so that 
$(x_i,y_j)=\kron i j$ and $g$ interchanges $x_i$ and $y_i$ for $i=1, \dots , r$ and fixes each $x_j, y_j$ for $j \ge r+1$.  
\end{lem}  
\pf  Let $W$ be the codimension 1 subspace of $[V,g]$ which contains all the singular vectors of $[V,g]$.  Take a basis $u_i$, $i=1,\dots , 2k$, of $[V,g]$ of nonsingular vectors. 
For $x \in [V,g]$, let $P(x):=\{v \in V|v(g-1)=x\}$, a coset of $[V,g]^\perp$.  
For all $x$, $P(x)$ contains singular vectors (see \reftt{annihsing}).  We therefore may take $x_1$ so that $x_1(g-1)=u_1$ and we define $y_1:=x_1^g$.  We may use induction on $span\{x_1,y_1\}^\perp$.  The only problem might be that we are unable to use \reftt{annihsing} at the last stage in case $r=\frac d 2$.  But then we use the fact that $V$ has plus type and the conclusion is forced.  
\eop

\begin{lem}\labtt{liftclass} (i) Suppose that $t$ is a clean upper involution of $\gd d$.  
Then the coset $t\rd d$ represents $s+1$ different conjugacy classes of involutions in $\gd d$, 
where $s$ is the number of orbits of $C_{\gd d}(t)$ on the cosets of 
$[R,t]$ in 
$C_{\rd d}(t)$ which contain involutions.   
We have $s=1$ if $k=\frac d 2$ and $s=2$ if $k<\frac d 2$.  This gives respectively one and two dirty classes of involutions in the coset.  

(ii) If $t$ is $\eta_{d,2k,\pm}$ (so is dirty and nonsplit), the coset 
$t\rd d$ represents one class of involutions if $k=\frac d 2$, and two  otherwise; all involutions in $t\rd d$ are dirty.     
\end{lem} 
\pf
Exercise.  
\eop

\begin{lem}\labtt{alldefectsoccur} 
(i) 
A defect $k$ involution in $\gd d/ \rd d \cong \Omega^+(2d, 2)$ is represented in $\brw d$ by an involution, specifically, by either 
a clean involution of defect $k$, or the dirty nonsplit involution $\eta_{d,2k,+}$, for a unique integer $k\le \frac d 2$.  
Furthermore, for any $d$ and positive   $k \le \frac d 2$, both cases occur and are mutually exclusive.  

(ii) An eigenlattice of $\eta_{2,2,+}$ has an orthogonal basis, of norms 2, 4.  
\end{lem}
\pf  
It is clear from a direct construction (or \reftt{alldefects}) and \reftt{liftclass} that both cases occur and that they are mutually exclusive.  
Since $\gd d$ contains a natural central product of $k$ natural $\brw 2 \cong W_{F_4}$ subgroups, it suffices to give a direct construction for the case $k=\frac d 2 =1$, which we now do.  
Notice that for $d=2$, $\brw 2\cong W_{F_4}$ contains two conjugacy classes of  reflection (upper and clean, of defect 1, representing the two classes when $k=\frac d 2$) and a nonsplit involution. 
Note that  the product of two reflections for orthogonal roots has trace 0, so is dirty.  There are two  orbits of $W_{F_4}$ on orthogonal pairs of roots, distinguished by root lengths, but the resulting products of two reflections represent only two classes: one class (for the pairs of equal length roots) and a second class for the case of unequal root lengths.  The latter gives the upper class.  For this case, we have an orthogonal set of vectors of norms 2 and 4 in a given eigenlattice, $M$, corresponding to orthogonal roots of different lengths.  
\eop

\bigbreak



\begin{thebibliography}{GRC99}






\bibitem[BW]{BW} E. S. Barnes  and G. E. Wall, Some extreme forms
defined in terms of abelian groups, JAMS 1 (1959), 47-63.   

\bibitem[BRW1]{BRW1} Beverly Bolt, T. G. Room and G. E. Wall, 
On the Clifford Collineations, Transform and Similarity Groups, I. 
Journal of the Australian Mathematical Society, 2, 1961, 60-79.  


\bibitem[BRW2]{BRW2} Beverly Bolt, T. G. Room and G. E. Wall, 
On the Clifford Collineations, Transform and Similarity Groups, II.
Journal of the Australian Mathematical Society, 1961, 80-96.  

\bibitem[B]{B} Beverly Bolt, T. G. Room and G. E. Wall, 
On the Clifford Collineations, Transform and Similarity Groups, III;
Generators and Relations,  Journal of the Australian Mathematical Society,
1961

\bibitem[Bour]{Bour} N. Bourbaki, \'Elements de Math\' ematique,
Groupes  et alg\`ebres de Lie, Chapitres 2 et 3, Diffusion C.C. L.
S., Paris, 1972.  


\bibitem[BE]{BE} Michel Brou\'e and Michel Enguehard, Une famille
infinie de formes quadratiques enti\`ere; leurs groupes d'automorphismes,
Ann. scient. \'Ec. Norm. Sup., $4^{eme}$ s\'erie, t. 6, 1973, 17-52. 

\bibitem[Car]{Car}  Roger Carter, Simple Groups of Lie Type,
Wiley-Interscience, London (1972).  




\bibitem [CR]{CR} Charles Curtis and Irving Reiner, Representation Theory 
of Groups and Associative Algebras, Interscience, 1962.  

\bibitem [Dieud]{Dieud} Jean Dieudonne, La G\'eom\'etrie des Groupes Classiques, Springer, Berlin Heidelberg New York, 1971.  
%




\bibitem[Gor]{Gor} 
 Daniel Gorenstein, Finite Groups, Harper and Row, New York, 1968.  

%

\bibitem[Gr72]{Gr72}  Robert L. Griess, Jr., 
 Automorphisms of extra special groups and nonvanishing degree 2 cohomology
(research announcement for [5]), in Finite Groups 1972:
Proceedings of the Gainesville Conference on Finite Groups, (T. Gagen, M.
P. Hale and E. E. Shult, eds.), North Holland Publishing Co., Amsterdam,
68-73, 1973.


\bibitem[GrDemp]{GrDemp} Robert L. Griess, Jr., On a subgroup of order 
$2^{15}|GL(5,2)|$  in $E_8 (C)$,  the Dempwolff
group and $Aut(D_8 \circ D_8 \circ D_8)$ , J. Algebra, 40, 1976, 271-279.


\bibitem[GrElAb]{GrElAb}   Robert L. Griess, Jr., Elementary abelian
subgroups of algebraic groups, Geometria Dedicata, {\bf 39}, 253-305, 1991.  


	
	\bibitem[GrEx]{GrEx}Robert L. Griess, Jr.,  Automorphisms of extra
	special groups and nonvanishing degree 2 cohomology, Pacific J. Math., 48,
	403-422, 1973.

\bibitem[GrNW]{GrNW}  
 Robert L. Griess, Jr., 
Sporadic groups, code loops and nonvanishing cohomology, J. Pure Appl.
Algebra, 44, 1987, 191-214.



\bibitem[GrMont]{GrMont}    Robert L. Griess, Jr., 
The monster and its nonassociative algebra, in Proceedings of the Montreal
Conference on Finite Groups, Contemporary Mathematics, 45, 121-157, 1985, 
American Mathematical Society, Providence, RI.

%

\bibitem[G12]{G12} Robert L. Griess, Jr.,   Twelve Sporadic Groups,
Springer Verlag, 1998.  

%
\bibitem[POE]{POE} Robert L. Griess, Jr, Pieces of Eight, Advances in
Mathematics, 148, 75-104 (1999).

\bibitem[PO$2^d$]{bwy} Robert L. Griess, Jr.,  Pieces of $2^d$: existence and uniqueness for
Barnes-Wall and Ypsilanti lattices.    56 pages.  To appear in Advances in Mathematics.  ; see also 


\bibitem[PO$2^d$corr]{bwycorr} Robert L. Griess, Jr., 
Corrections and additions to `` Pieces of $2^d$: existence and uniqueness for
Barnes-Wall and Ypsilanti lattices. '', which will be posted on the author's web site and on arxiv.  



\bibitem[Hup]{Hup} Bertram Huppert,  Endliche Gruppen I, Springer Verlag,  
Berlin, 1968. 



\bibitem[MS]{MS}
Jesse MacWilliams and Neal Sloane, The Theory of Error Correcting Codes,
North-Holland, 1977.

\bibitem[McL]{McL} Jack E. McLaughlin, Some subgroups generated by transvections, Arch. Math. 18 (1969), 108-115.  

\bibitem[MH]{MH} Milnor and Husemoller,  Symmetric Bilinear Forms,
Ergebnisse der Mathematick und Ihrer Grenzgebiete, Band 73, Springer
Verlag, New York, 1973.  

\bibitem[Poll]{Poll} Harriet Pollatsek, Cohomology groups of some linear groups over fields of characteristic 2, Illinois Journal of Mathematics, 15 (1971) 393-417.  

\bibitem[Se]{Se} Jean-Pierre Serre, A Course in Arithmetic, Springer
Verlag, Graduate Texts in Mathematics 7, 1973.  

%

\end{thebibliography}
\end{document}